\documentclass[11pt]{amsart}
\usepackage{amsmath,amssymb,latexsym,soul,cite,amsthm,color,enumitem,graphicx,tikz, mathtools,microtype}
\usepackage[export]{adjustbox}

\usepackage{amsfonts,amsmath,amsthm,latexsym,amssymb,marvosym,esint,xfrac,bbm,amssymb}

\usepackage{cite}

\usepackage{graphics, epsfig}
\tikzstyle{mybox} = [draw=black, very thick, rectangle, rounded corners, inner ysep=5pt, inner xsep=5pt]

\usepackage[colorlinks=true, urlcolor=blue, citecolor=red, linkcolor=blue]{hyperref}

\usepackage[english]{babel}
\usepackage[left=2.5cm,right=2.5cm,top=2.7cm,bottom=2.7cm]{geometry}

\numberwithin{equation}{section}

\newcommand{\osc}{\operatornamewithlimits{osc}}

\newcommand{\dist}{\operatorname{dist}}

\newtheorem{theorem}{Theorem}[section]
\theoremstyle{plain}
\newtheorem{lemma}[theorem]{Lemma}
\theoremstyle{plain}

\theoremstyle{plain}
\newtheorem{corollary}[theorem]{Corollary}
\newtheorem{definition}[theorem]{Definition}
\theoremstyle{definition}
\newtheorem{remark}[theorem]{Remark}

\newcommand{\N}{{\mathbb N}}

\newcommand{\R}{{\mathbb R}}

\newcommand{\e}{\epsilon}

\newcommand{\beq}{\begin{equation}}
\newcommand{\eeq}{\end{equation}}

\renewcommand{\ge}{\geqslant}

\def\Xint#1{\mathchoice
    {\XXint\displaystyle\textstyle{#1}}%
    {\XXint\textstyle\scriptstyle{#1}}%
    {\XXint\scriptstyle\scriptscriptstyle{#1}}%
    {\XXint\scriptscriptstyle\scriptscriptstyle{#1}}%
    \!\int}
\def\XXint#1#2#3{\setbox0=\hbox{$#1{#2#3}{\int}$}
    \vcenter{\hbox{$#2#3$}}\kern-0.5\wd0}
\def\bint{\Xint-}
\def\bint{\Xint-}
\def\dashint{\Xint{\raise4pt\hbox to7pt{\hrulefill}}}
\def\dashiint{\bint\kern-0.15cm\bint}

\def\tr{(u-k)_{-}}
\def\A{\mathbb{A}}
\def\Q{\mathcal{Q}}
\def\K{\mathcal{K}}
\def\dive{\mathrm{div}}
\def\d{\mathrm{d}}
 \def\B{\mathcal{B}}
\def\dim{\mathrm{dim}}
\def\dist{\mathrm{dist}}
\def\H{\mathcal{H}}
\def\I{\mathcal{I}}
\def\L{\mathcal{L}}

\def\P{\mathcal{P}}
\def\T{\mathcal{T}}
\def\S{\mathcal{S}}

\title[Wiener's criterion for degenerate double-phase diffusion] {Fine boundary continuity for degenerate double-phase diffusion}
\author[Ciani, Henriques, Skrypnik]{Simone Ciani {$\&$} Eurica Henriques {$\&$} Igor I. Skrypnik}
\address{Department of Mathematics of the University of Bologna, Piazza Porta San Donato, 5, 40126 Bologna, Italy}

\email{simone.ciani3@unibo.it}

\address{Centro de Matem\'atica, Universidade do Minho - Polo CMAT-UTAD Departamento de Matem\'atica Universidade de Tr\'as-os-Montes e Alto Douro, Vila Real, Portugal}

\email{eurica@utad.pt}

\address{Institute of Applied Mathematics and Mechanics, National Academy of Sciences of Ukraine, Gen. Batiouk Str. 19, 84116 Sloviansk, Ukraine}

\email{ihor.skrypnik@gmail.com}

\begin{document}
\begin{abstract} 
We study the boundary behavior of solutions to parabolic double-phase equations through the celebrated Wiener's sufficiency criterion. The analysis is conducted for cylindrical domains and the regularity up to the lateral boundary is shown in terms of either its $p$ or $q$ capacity, depending on whether the phase vanishes at the boundary or not. Eventually we obtain a fine boundary estimate that, when considering uniform geometric conditions as density or fatness, leads us to the boundary H\"older continuity of solutions. In particular, the double-phase elicits new questions on the definition of an adapted capacity.

\vskip0.2cm \noindent

\end{abstract}

\maketitle

 \begin{center}
		\begin{minipage}{9cm}
			\small
   \tableofcontents
		\end{minipage}
	\end{center}



\def\tr{(u-k)_{-}}
\def\R{\mathbb{R}}
\def\A{\mathbb{A}}
\def\Q{\mathcal{Q}}
\def\b{\mathbb{B}}
\def\K{\mathcal{K}}
\def\k{\mathbb{K}}
\def\E{\mathcal{E}}
\def\q{\mathbb{Q}}

\def\N{\mathbb{N}}
\def\dive{\mathrm{div}}
\def\d{\mathrm{d}}
 \def\B{\mathcal{B}}
\def\dim{\mathrm{dim}}
\def\dist{\mathrm{dist}}
\def\H{\mathcal{H}}
\def\I{\mathcal{I}}
\def\L{\mathcal{L}}
\def\C{\mathcal{C}}
\def\T{\mathcal{T}}
\def\S{\mathcal{S}}
\def\P{\mathcal{P}}
\def\e{\epsilon}


\addtocontents{toc}{\protect\setcounter{tocdepth}{1}}

\section{Introduction and main results}\label{Introduction}
\noindent Let $u(x,t)$ describe the flow, in a space configuration $x \in \Omega \subset \R^N$ and at a time $t \in (0,T)$, of the velocity of a non-Newtonian fluid which changes the power-law of its stress tensor according to a dramatic switch of the energy density. This change is specified by a law $a(x,t)$, that can be catered by an electromagnetic field or a mechanical device that suddenly obstructs the flow. Some of the fluids just described are addressed as electro-rheological and are of promising technological interest (see for instance \cite{RajRuz}, \cite{RajWin} or the book \cite{Ruz}); their special feature being a heavy change of viscosity in a very short time. As a guiding example, in \cite{ApuBilFuc} and \cite{BilFuc1}, the authors consider the stationary flow of a generalized non-Newtonian fluid, modeled after an anisotropic dissipative potential $\Phi(z)=|z|^p+a(x,t) |z|^q$, whose energy is trapped between two power laws. Here we are interested in this description, as opposed to a slower change of rate that happens when $\Phi(z)= |z|^{p(x)}$ and $p$ is a log-continuous function. In the former case, the regularity of the solution, if any, is expected to follow a rule dictated by $a(x,t)$ itself, which we call {\it the phase}. In the present work we propose an analysis of the boundary behavior of solutions to equations that embody these features, and whose prototype is referred to as the parabolic double-phase equation, given by
\begin{equation} \label{proto} \partial_t u - \dive \bigg( \left(|\nabla u |^{p-2}+a(x,t) |\nabla u|^{q-2}\right) \nabla u \bigg)=0 \quad \text{in} \quad \Omega_T=\Omega \times (0,T] ,\end{equation}
\noindent for $\Omega \subset \R^N$ open and bounded. Given a continuous initial datum $f$ prescribed on the parabolic boundary of $\Omega_T$, we address the question of whether solutions $u$ to the parabolic double-phase equation \eqref{proto} reach $f$ in a continuous fashion; and in such case, when the datum is more regular (for instance, H\"older continuous), we would like to describe how fast this happens. Our answer to this question is presented in Theorem \ref{boundary-estimate}. In particular, we find that the trade-off between the geometry of the lateral boundary $\partial \Omega \times (0,T)$ in terms of the elliptic $p$ or $q$ capacity of $\partial \Omega$ and the behaviour of the phase at these points determines the desired rate.

\subsection{Origins and Framing of the Topic}
\noindent In recent years, the stationary version of equation \eqref{proto} has received a great attention, especially with regard to the regularity theory; we refer, for instance, to the surveys \cite{Mar-surv}, \cite{MinRad}, \cite{Pap}, and the extensive lists of references therein. While the local boundedness of solutions was already studied in the 70's in \cite{Kol} and \cite{Kor}, the non-standard behavior was faced by Zikhov in \cite{Zik} in the context of averaging of variational problems and the first pioneering analysis of the regularity of the gradient appeared in \cite{Mar1}, \cite{Ura} (see also \cite{BraBou} for a full-anisotropic version). In parallel, a fruitful theory of adapted and generalized energy spaces has seen the light, as Orlicz and Musielak-Orlicz spaces: here we refer, for instance, to the practical survey \cite{Chle}.

\noindent Equation \eqref{proto} belongs to a wider class of equations exhibiting the so-called ($p,q$)-growth, that for its mathematical challenges together with its numerous applications draw a considerable attention for several decades. Regarding the stationary point of view, a non-exhaustive list of contributions is  
\cite{AceMinSer, Alk, AlkSur, ArrHue, BarColMin1, BarColMin2, BarColMin3, BenHarHasKar, ColMin1, ColMin2, ColMin3, CupMarMas, CupMarMas2, HadSkrVoi, HarKinLuk, HarHasLee, HarKuuLukMarPar, Ok, RagTac,  SkrVoi1} to which we refer for results, references, historical notes and extensive survey of regularity issues, being the literature so wide that it results complicated to track every result in this direction. \vskip0.1cm  
\noindent On the other hand, the regularity theory for evolutionary double-phase equations has received less attention, most probably because of the merging of the difficulties inherent to the double-phase with the ones of the non-homogeneity of the operator caused by the parabolic term. A study of the local $L^{\infty}$ norm of the gradient has been brought on in \cite{BarLin}, \cite{BoeDuzMar}, \cite{DieSchSch} and \cite{Sin2}. Refined quantitative gradient bounds have been addressed in \cite{Cri}, while higher differentiability of the gradient has been investigated in \cite{GiaPasSch}. 

\vspace{0.2cm}

\noindent Our interest specifies towards equations with measurable and bounded coefficients, in the framework of a fine boundary estimate that is irrespective of the higher-order regularity. Within this perspective, the continuity and H\"older continuity for parabolic equations with Orlicz growth (generalizing \eqref{proto}) has been studied in \cite{BarLin}, \cite{HwaLie1}, \cite{HwaLie2}, \cite{Sin}, \cite{SkrVoi1} and \cite{SkrVoi2}; while in \cite{BurSkr}, \cite{SavSkrYev1} and \cite{SavSkrYev2} the authors proved suitable versions of the Harnack inequality (see also \cite{Sur}, \cite{Wan} for the variable exponent case).

\subsection{Fine Boundary Regularity}  A sufficient condition for the regularity of a boundary point for the prototype $p$-Laplacian elliptic equation has been known since the famous paper of Maz'ya \cite{Maz}, and is named after Wiener, who studied the Dirichlet problem for the linear case from the potential point of view (see \cite{Wiener1}, \cite{Wiener2}). Later, Gariepy and Ziemer in \cite{GarZie} generalized this criterion to the case of quasi-linear elliptic equations. Roughly speaking this sufficiency condition is the following: picking $x_o \in \partial{\Omega}$ and defining for $p>1$, $r>0$ the number
\begin{equation}\label{delta}
    \delta_p(r)= \bigg( \frac{\C_p(\overline{B_r(x_o)}\setminus \Omega; B_{2r}(x_o))}{\C_p(\overline{B_{r}(x_o)}; B_{2r}(x_o))} \bigg)^{\frac{1}{p-1}},
\end{equation} \noindent where $\C_p(K,B)$ is the elliptic variational $p$-capacity of the condenser $(K,B)$ (see \eqref{capacity} for details), then, weak solutions of quasi-linear elliptic equations of $p$-Laplacian type are continuous up to the point $x_o$ if 
\begin{equation}\label{eq1.3}
\int_{0}^1 \delta_p(r)\frac{dr}{r}=\infty\, .
\end{equation}
\noindent We will refer, here and in the sequel, to the books \cite{HeiKipMar} and\cite{MalZie} for an account of capacity methods for the fine boundary regularity in the context of elliptic $p$-Laplacian type equations.\vskip0.1cm  
\noindent The problem of fine boundary regularity for the diffusive $p$-Laplacean equation is much more recent. Continuity up to the boundary with monotonicity conditions was proved in \cite{Skr1}, \cite{Skr} under the condition \eqref{eq1.3}. This result was generalized in \cite{GiaLia} for more general parabolic evolution equations by using a weak Harnack inequality (see also \cite{GiaLia1}, \cite{GiaLiaLuk} for the singular super/sub-critical cases). \vskip0.1cm \noindent Finally, a sufficiency criterion of Wiener-type for parabolic equations with non-standard growth conditions is, up to our knowledge, a novelty. The present work is therefore a first step on the understanding of boundary regularity for non-standard parabolic operators.

\subsection{Setting of the Problem}
\noindent Let us denote $B_\rho(y)$ the ball in $\R^N$ of center $y$ and radius $\rho$, and let $\Omega$ be a bounded domain in $\mathbb{R}^{N}$. For $T>0$, we consider $\Omega_{T}:= \Omega \times (0, T]$ the cylinder with base $\Omega$ and length $T$, and we denote by $S_{T}:=\partial \Omega\times (0,T]$ its lateral boundary. We consider equations 
\begin{equation}\label{eq1.1}
\partial_t u-\textrm{div}\, \mathbb{A}(x, t, \nabla u)=0, \quad\text{weakly in} \ \Omega_{T},
\end{equation}
where we assume that the function $\mathbb{A}:\Omega_{T}\times \mathbb{R}^{N} \rightarrow \mathbb{R}^{N}$ is Caratheodory, i.e. $\mathbb{A}(\cdot, \cdot, \xi)$ is Lebesgue measurable for all $\xi \in \mathbb{R}^{N}$,
and $\mathbb{A}(x, t, \cdot)$ is continuous for almost all $(x, t)\in \Omega_{T}$; and that $\mathbb{A}$ satisfies the following structure conditions
\begin{equation}\label{eq1.2}
\begin{aligned}
\mathbb{A}(x, t, \xi)\cdot \xi &\geqslant  C_{1}\,\big( |\xi|^{p} + a(x,t)|\xi|^{q}\big)=:C_{1}\,\varphi(x,t, |\xi|), \quad 2< p < q,
   \\
\ |\mathbb{A}(x, t, \xi)| &\leqslant  C_{2}\big( |\xi|^{p-1} + a(x,t)|\xi|^{q-1}\big)=: C_{2}\,{\varphi(x,t,|\xi|)}/{|\xi|},
\end{aligned}
\end{equation}
for $C_{1}, C_{2}$ given positive constants, that we will refer to as structural data. In addition, we assume that the function $a(x,t): \R^{N+1}\rightarrow [0,\, \infty)$ is everywhere defined and non-negative. We assume $a(x,t)$ to be locally H\"older continuous around $S_T$: for any $(x_{o}, t_{o})\in S_{T}$, we assume that there exist positive numbers $R_o,A_o$ such that, for any $0<r<R_o$, the following inequality holds true,
\begin{equation} \label{A}
\osc_{Q_{r,r^{2}}(x_{o},t_{o})} a(x,t)\leq A_o \ r^{q-p},
\end{equation} \noindent being $Q_{r,r^{2}}(x_{o},t_{o})= B_r(x_o) \times (t_o-r^2, t_o+r^2)$.

\noindent As our estimates are local in nature, the constants $R_o$ and $ A_o$ will also be referred to as structural constants. Thence, we are concerned with the boundary behaviour of solutions to the Cauchy-Dirichlet problem
\begin{equation}\label{C-pb}
    \begin{cases} \partial_t u - \dive \,  \A(x,t,\nabla u)=0, & \text{weakly in} \ \ \Omega_T,\\
    u(x,t)= f(x,t), & \text{on} \,  S_T,\\
   u(x, 0)= f(x,0), & \text{attained in} \, \, L^2_{loc}(\Omega),
\end{cases} \end{equation}\noindent where $\mathbb{A}$ obeys to \eqref{eq1.2}-\eqref{A} above for $2<p<q$, and 
\[f \in L^q(0,T; W^{1,q}(\Omega))\cap C(\overline{\Omega_T}).\] 
\noindent The boundary datum $f$ is taken in the weak sense, i.e. $(u-f)(\cdot, t) \in W^{1,q}_o(\Omega)$ for almost every time $t \in (0,T]$. As typical of parabolic equations, what happens in the future is determined entirely from the past: this motivates the omission of a prescription of the boundary datum at $\Omega \times \{T\}$. In agreement with this principle, for our local estimates we will work with backward parabolic cylinders (See Section \ref{Notation} for more details).\newline  
\noindent The well-posedness of this problem has been addressed in \cite{BoeDuzMar}, \cite{Sin} and very recently in \cite{AroShm}, with slightly different notions of solutions. We refer to Section \ref{Preliminaries} below for the details of our definitions.\vskip0.1cm 

\noindent Finally, another important topic concerns global boundedness of solutions, for which there seems not to be a complete picture in the parabolic case for equations such as \eqref{eq1.1}. In general and within an elliptic context, for this generality of choice of exponents $q>p>2$, local weak solutions to stationary equations with $(p,q)$ growth as \eqref{eq1.1} above are not meant to be locally bounded, as the two pioneering counter-examples \cite{Gia-cex}, \cite{Mar-cex} show. Nonetheless, these two examples are fully anisotropic, meaning with this that the energy is not a function of the modulus of the gradient, but just of its components. For general non-standard parabolic equations, global boundedness is shown in Theorem 3 of \cite{MinXit} for fully anisotropic parabolic equations; see also \cite{CiaVesVes} for refined local bounds. The condition given may not be sharp in the case of equation \eqref{eq1.1}, as it is unrelated to the degree of H\"older continuity of $a(x,t)$ (see for instance\cite{Sin}), or its $L^{\infty}$ norm. For this reason, in what follows we consider solutions that are globally bounded in $\Omega_T$, thereby admitting a wider set of solutions.

\subsection{Main Result and Applications} \noindent In order to formulate our boundary estimate, we briefly recall the definition of capacity at hands. Let $s \in (1,N]$, $B \subset \R^N$ be an open set and $K\subset B$ be a compact set. We denote by $\C_s(K; B)$ the Newtonian (or variational) capacity of the condenser $(K;B)$ and defined as
\begin{equation}\label{capacity}
\C_s (K;B)= \inf \bigg{\{}\|\nabla f\|_{L^s(B)}^s\,:\, f \in C_o^{\infty}(B),\, \, f \ge 1\, \, \text{on}\, \, K  \bigg{\}}.
\end{equation} \noindent This introduced version of capacity pertains to domains of $\R^N$ and it extends, within its elliptic fashion, spontaneously to cylinders $Q=B \times (t_1,t_2)$ in $\R^{N+1}$. Let $\tilde{K} \subset Q$ be a compact subset of such cylinder, and if we define
\[\mathcal{C}_s(\tilde{K},Q) = \inf \bigg{\{}\|\nabla f\|_{L^s(Q)}^s\, : \, \, f \in C_o^{\infty}(Q), \, \, f \ge 1 \, \, \text{on} \, \, \tilde{K} \bigg{\}}, \quad \text{then} \quad \mathcal{C}_s(\tilde{K},Q) = \int_{t_1}^{t_2} \mathcal{C}_s (\tilde{K}_{\tau},B)\, d \tau,\]
\noindent being $\tilde{K}_\tau= K \times \{\tau\}$. The proof of this last equality can be found in \cite{BirMos}, while other notions of parabolic capacity are investigated in various other circumstances, see for instance \cite{AveKuuPar} and \cite{Zie}. \vskip0.1cm  

\noindent With this definition and \eqref{delta}, our main result reads as follows.
\begin{theorem}\label{boundary-estimate} Let $(x_o,t_o) \in S_T$ and let $u$ be a bounded, weak solution to the Cauchy-Dirichlet problem \eqref{C-pb}. Depending on the point $(x_o,t_o)$, we assume that either
\begin{equation}\label{eq1.6}
\int_0^1 \delta_p(r) \frac{dr}{r}=\infty ,  \quad  \text{if} \quad a(x_o,t_o)=0,
\end{equation}
or
\begin{equation} \label{eq1.7}
\int_0^1 \delta_q(r) \frac{dr}{r}=\infty ,    \quad \text{if} \quad a(x_o,t_o)>0.
\end{equation} \noindent Then, in each case respectively, there exist $\{\rho_0(p),\eta_0(p)\}$, $\{  \rho_0(q), \eta_0(q)\}$ couples of positive numbers depending only on the data and conditions \eqref{eq1.6}-\eqref{eq1.7}, and positive constants $\gamma, \hat{\gamma}, \gamma^*$ depending only on the data, such that, defining 
\[Q_0(p)=B_{\rho_0(p)}(x_o)\times (t_o-\eta_0(p), \, t_o],\]
\[Q_0(q)=B_{\rho_0(q)}(x_o)\times (t_o-\eta_0(q), \, t_o], \quad \text{and} \quad \omega_0= \osc_{\Omega_T} u,\]
\noindent the following inclusions 
\[Q_{\rho}(\omega_0,p)= B_{\rho}(x_o) \times (t_o-\gamma^* \rho^p \omega_0^{2-p},\,  t_o] \subset Q_0(p),  \]
\[Q_{\rho}(\omega_0,q)= B_{\rho}(x_o) \times (t_o-\gamma^* \rho^q \omega_0^{2-q},\,  t_o] \subset Q_0(q), \]
\noindent hold true for $\rho=\rho_0(p), \rho_0(q)$, and for all $0<\rho<\rho_0(p)$ we have the estimate
\begin{equation*}\label{resultp}
     \osc_{Q_{\rho}(\omega_0,p)\cap \Omega_T} u \leq \omega_0 \exp \bigg{\{} -\frac{1}{\gamma} \int_{\rho}^{\rho_0(p)} \delta_p (s) \frac{ds}{s} \bigg{\}}+ \osc_{Q_{0}(p) \cap S_T} f + \hat{\gamma} [\rho_0(p)]^{\frac{\epsilon}{p-2}}, \quad\text{if} \quad a(x_o,t_o)=0,
\end{equation*} 
\noindent while for all $0<\rho< \rho_0(q)$ we have 
\begin{equation*}\label{resultq}
     \osc_{Q_{\rho}(\omega_0,q)\cap \Omega_T} u \leq \omega_0 \exp \bigg{\{} -\frac{1}{\gamma} \int_{\rho}^{\rho_0(q)} \delta_q (s) \frac{ds}{s} \bigg{\}}+ \osc_{Q_{0}(q) \cap S_T} f + \hat{\gamma} [\rho_0(q)]^{\frac{\epsilon}{q-2}}, \quad\text{if} \quad a(x_o,t_o)>0.
\end{equation*} 
\end{theorem}

\noindent We observe that the geometric construction is dependent on the assumptions \eqref{eq1.6}-\eqref{eq1.7}, differently from the isotropic singular case (see for instance \cite{GiaLiaLuk}). Nonetheless, even if Theorem \ref{boundary-estimate} is stated for the Cauchy-Dirichlet problem, as soon as a lateral boundary datum is concerned, the oscillation estimates above are of local nature. \vskip0.1cm \noindent 
In this framework, it is a simple consequence  that a Wiener-type test is a sufficient condition for a point $(x_o,t_o)\in S_T$ to be a regular point to the parabolic double-phase operator \eqref{eq1.1}-\eqref{eq1.2}-\eqref{A}.\vskip0.1cm \noindent We recall that a lateral boundary point $(x_{o}, t_{o})\in S_{T}$ is said to be regular to  \eqref{eq1.1}-\eqref{eq1.2}-\eqref{A} if, for any weak solution $u$ to equation \eqref{eq1.1}, satisfying 
\begin{equation}\label{bd-datum}
(u(x,t)-f(x,t))\in V^{2,q}_{o}(\Omega_{T}), 
\end{equation} 
with any $f(x,t) \in C(\overline{\Omega_{T}})$, the limit
\begin{equation*}
\lim\limits_{\Omega_{T} \ni(x,t)\rightarrow(x_{o},t_{o})}u(x,t)=f(x_{o}, t_{o})
\end{equation*} \noindent 
is attained. Here and in what follows, we denote with $V_o^{2,q}$ the parabolic space \[
  V^{2,q}_o(\Omega_{T})= C(0, T; L^{2}(\Omega))\cap L^{q}(0, T; W_{o}^{1,q}(\Omega)),\] and the attainment of the datum  \eqref{bd-datum} is understood weakly. The geometric conditions \eqref{eq1.6}-\eqref{eq1.7} are also common in the literature when referring to the set $\R^N\setminus{\Omega}$ as ($p$ or) $q$-thick at $x_o$ (e.g. \cite{HeiKipMar}). 

\noindent 
\begin{corollary} \label{th1.1}
Let $u$ be a bounded, weak solution to equation \eqref{eq1.1}-\eqref{eq1.2}, and let \eqref{A} be satisfied in $(x_o,t_o)\in S_T$. If moreover 
\vskip0.2cm \noindent 

\begin{itemize}
    \item $a(x_o,t_o)=0$, then \eqref{eq1.6} is a sufficient condition for $(x_{o}, t_{o})$ to be regular to \eqref{eq1.1}-\eqref{eq1.2}-\eqref{A};
\end{itemize}

\noindent otherwise, if    
\begin{itemize}
    \item   $a(x_{o}, t_{o}) >0$, 
then \eqref{eq1.7} is a sufficient condition for $(x_{o}, t_{o})$ to be regular to \eqref{eq1.1}-\eqref{eq1.2}-\eqref{A}.
\end{itemize} 
\end{corollary}
\noindent Classically, in the case $p=q=2$, when at the point $x_o\in \Omega$ further requirements are satisfied, as the logarithmic Wiener condition (see \cite{Bir}), the solutions attain a H\"older continuous datum in a H\"older continuous fashion. For ease of exposition here we ask $\Omega$ to enjoy a uniform geometrical property; which is ensured, for instance, by the classic corkscrew condition (see \cite{HeiKipMar} Thm 6.31). We briefly recall it here below.\vskip0.1cm 

\noindent Let $X \subset \R^N$ be a closed set, $Y \subseteq X$ and $s \in (1,N]$. We recall that the set $X$ is uniformly $s$-fat in $Y$ if there exist positive constants $\lambda_s, R_s$ such that, for all $y \in Y$ and $0<\rho<R_s$,
\[\C_s( \overline{B_\rho(y)} \cap X; B_{2\rho}(y)) \ge \lambda \rho^{N-s}.\]
\noindent When $X=Y$ we just say that $X$ is uniformly $s$-fat. With this definition at hand, we can present a notion of fatness that suits the double-phase problem.

\begin{definition}
    Given a continuous function $a:\R^{N+1}\rightarrow [0,\,\infty)$, we say that a closed set $X \subset \R^N$ is uniformly $(p,q)$-fat with phase $a(x,t)$ if $X$ is uniformly $p$-fat at those points $x_o \in \partial X$ such that $a(x_o,t_o)=0$ for some $t_o \in \R$, and it is uniformly $q$-fat at those points $x_o \in \partial X$ such that $a(x_o,t)>0$ for all $t \in \R$.
\end{definition}
\begin{remark} We observe that in the above definition if $X$ is uniformly  $p$-fat and the function $a$ vanishes on $\partial X$ for all times, then trivially $X$ is uniformly $(p,q)$-fat with phase $a(x,t)$ with any $q$. Moreover, when $q\ge p$, a uniformly $p$-fat set is also a uniformly $q$-fat set, by a simple application of H\"older's inequality. Hence the introduced definition is weaker than the usual $p$-fatness. The property of a set of being uniformly $p$-fat is an open-end condition (see for instance \cite{Lewis}) and it is equivalent to a point-wise Hardy inequality (see \cite{KorLehTuo}). The definition of fatness obliges $q<N$: in the cases where $p>N$ condition \eqref{eq1.6} is satisfied and when $q>N$ condition \eqref{eq1.7} is satisfied, because in such cases the capacities of point and ball are comparable with a uniform constant. This remark further implies that, when $p<N<q$, if for these times $t\in \R$ such that the set $\{X\times \{t\}\}\cap a(\cdot, t)^{-1}(\{0\})$ is not empty, it is also uniformly $p$-fat, then $X$ is uniformly $(p,q)$-fat with phase $a(x,t)$. 
\end{remark}

\noindent Finally, when the complement of $\Omega$ is uniformly $p$-fat, then the integral \eqref{eq1.6} diverges at every boundary point $x_o \in \partial \Omega$, and this leads us to the following corollary of Theorem \ref{boundary-estimate}.

\begin{corollary} \label{HC}
Let $u$ be a bounded, weak solution to \eqref{eq1.1}-\eqref{eq1.2}-\eqref{A}, with an H\"older continuous boundary datum $f \in C^{0,\alpha}(\overline{\Omega_T})$. Suppose furthermore that $\R^N\setminus \Omega$ is $(p,q)$-fat with phase $a(x,t)$. Then, the solution $u$ is H\"older continuous up to $S_T$.\end{corollary}

\noindent Classically (for instance in \cite{Campa}, \cite{Lady}, \cite{Trudy} and for the parabolic $p$-Laplacean in \cite{DiB}) the H\"older continuity up to the boundary was obtained for domains $\Omega \subset \subset \R^N$ satisfying the density condition
\begin{equation}\label{density}
    \exists \alpha, R_D>0\, : \, \forall x_o \in \partial \Omega \quad \forall \, 0<\rho<R_D\, \quad |\Omega \cap B_{\rho}(x_o)|\leq (1-\alpha) |B_\rho|.
\end{equation} \noindent By a simple application of the definition of the $s$-capacity together with the Poincar\'e inequality, condition \eqref{density} implies that $\R^N\setminus \Omega$ is uniformly $s$-fat; however, the converse statement is not true in general, as already seen by the case of points when $s>N$, or by the fact that sets of zero $s$-capacity do not separate the space $\R^N$. Nonetheless, when dealing with a global problem and for the purpose of precise integral estimates, these two conditions meet when a zero-extension is available; see for instance \cite{BCE}, Prop. 5.9 in the context of Campanato theory. Finally, we refer to Corollary 11.25 of \cite{BjoBjo} for more geometrical notions implying boundary regularity: among these examples, the $p$-fatness of the complement is the weakest assumption. 

\subsection*{Structure of the paper} In Section \ref{Notation}, we collect the notation used in the overall paper. Then, in Section \ref{Preliminaries}, we define local weak solutions and we describe various Lemmata concerning Energy (Caccioppoli) estimates, a measure-theoretical maximum principle, negative-powers Energy estimates, a Reverse H\"older's inequality and finally the weak Harnack inequality for nonnegative local weak supersolutions to \eqref{eq1.1}-\eqref{eq1.2}-\eqref{eq1.3}. In Section \ref{sect3}, we draw the geometric setting of the proof and we use the results of Section \ref{Preliminaries} to prove a reduction of oscillation of the solution near the boundary by means of the capacity of $\partial \Omega$ at the point considered. Finally in Section \ref{conclusion}, we prove the main result, Theorem \ref{boundary-estimate}, and in Section \ref{Appendix}, we collect the proof of the Energy Estimates of Section \ref{Preliminaries}, in order to leave space in the main text to what is really new. 

\section{Notation} \label{Notation}
{\small 
 \begin{itemize}
 \item {\it Constants dependency.} \newline \noindent We refer to the parameters $N$, $p$, $q$, $C_{1}$, $C_{2}$, $A_o$ and $M:=\sup\limits_{\Omega_{T}}|u|$ as our structural data, and we say that a constant $\gamma$ depends only on the data if it can be quantitatively determined a priori only in terms of the above quantities. The generic constant $\gamma$ may change from line to line. 
\vskip0.1cm \noindent 
     \item  {\it Geometry.} \newline \noindent We denote by $O$ the origin in $\R^N$. Let $r, \eta>0$. We denote with $B_r(x)$ the ball of radius $r$ centered in $x \in \R^N$.   Then we write 
     \begin{equation*} 
     \begin{cases} Q^{+}_{r,\eta}(\bar{x},\bar{t})= B_r(\bar{x}) \times (\bar{t}, \, \bar{t}+\eta),\\
     Q^{-}_{r,\eta}(\bar{x},\bar{t})= B_r(\bar{x}) \times (\bar{t}-\eta, \, \bar{t}),\\
     Q_{r,\eta}(\bar{x},\bar{t})= B_r(\bar{x}) \times (\bar{t}-\eta, \, \bar{t}+\eta), \end{cases} \end{equation*}
respectively, for the forward, backward and full cylinders centered at $(\bar{x},\bar{t})$ of radius $r$ and length $\eta$ (or $2\eta$). When writing
\[Q_{r}^{\pm}= Q_{r,r^2}^{\pm},\] we denote the cylinder centered at $O$ and whose time interval has length $r^2$; 
being \[Q_r= Q_{r,r^{2}} = Q_r^-\cup Q_r^+ \ .\]
\vskip0.1cm\noindent

\item {\it Levels.} \newline \noindent 
For any level $k \in \R$, $(\bar{x}, \bar{t}) \in \Omega_T$, $r,\eta$ as before such that the inclusion 
$ Q^{+}_{r,\eta}(\bar{x},\bar{t}) \subset \Omega_T$
is satisfied, we denote by:
 \[A^{-}_{k,r,\eta}= Q^{+}_{r,\eta}(\bar{x},\bar{t})\cap \big\{u \leqslant k\big\} \]
    the sub-level sets of $u$ in $Q^{+}_{r,\eta}(\bar{x},\bar{t})$  and by    
\begin{equation*} \label{phi}  \varphi^{\pm}_{Q^{+}_{r,\eta}(\bar{x},\bar{\eta})}\bigg(\frac{k}{r}\bigg)=\bigg(\frac{k}{r}\bigg)^{p}+ a^{\pm}_{Q^{+}_{r,\eta}(\bar{x},\bar{t})}\bigg(\frac{k}{r}\bigg)^{q},\end{equation*} 
\noindent where $a: Q \subset \R^{N+1} \rightarrow \R_0^+=[0,\, \infty)$,  $\displaystyle{a^{+}_{Q}=\max_{Q} a}$ and $ \displaystyle{a^{-}_{Q}=\min_{Q}  a }$   .

 \end{itemize}
 }


\section{Preliminaries}\label{Preliminaries}
\subsection{Definition of solution} \label{definition}We say that a function \[u \in V^{2,q}_{loc}(\Omega_{T}):= C_{\textrm{loc}}(0, T; L^{2}_{\textrm{loc}}(\Omega))\cap L_{\textrm{loc}}^{q}(0, T; W_{\textrm{loc}}^{1,q}(\Omega)),\]
is a local weak super(sub)-solution to \eqref{eq1.1} if for any compact set $E \subset \Omega$
and every sub-interval $[t_{1}, t_{2}]\subset (0, T]$ there holds
\begin{equation}\label{eq1.4}
\int_{E}u \zeta\, dx \bigg|^{t_{2}}_{t_{1}} + \int^{t_{2}}\limits_{t_{1}}\int_{E}
\{-u\partial_{\tau}\zeta+ \mathbb{A}(x, \tau, \nabla u) \nabla \zeta\}\, dx d\tau \geqslant\, 0,\,\qquad \quad(\leqslant 0),
\end{equation}
 for any nonnegative test function $\zeta \in W_{\textrm{loc}}^{1,2}(0, T; L^{2}(E))\cap L_{\textrm{loc}}^{q}(0, T; W_{o}^{1,q}(E))$.\newline 
 \noindent A function \[u \in  C(0, T; L^{2}(\Omega))\cap L^{q}(0, T; W^{1,q}(\Omega)),\] such that 
 \[(u-f)\in W_o^{1,q}(\Omega) \quad \text{for a.e.} \quad t \in (0,T],\]
 is a weak super(sub)-solution to the Cauchy-Dirichlet problem \eqref{C-pb}, if for all $t \subset (0, T]$ it satisfies
\begin{equation}\label{defsol-CP}
\int_{\Omega}u \zeta(x,t)\, dx + \iint_{\Omega_T}
\{-u\partial_{\tau}\zeta+ \mathbb{A}(x, \tau, \nabla u) \nabla \zeta\}\, dx d\tau \geqslant\,  \int_{\Omega} f\zeta(x,0)\, dx,\,\qquad \quad(\leqslant 0),
\end{equation}
 for any nonnegative test function $\zeta \in W^{1,2}(0, T; L^{2}(\Omega))\cap L^{q}(0, T; W_{o}^{1,q}(\Omega))$.

 \vskip0.2cm 
 
 \noindent To the aim of our computations, it is technically convenient to have a formulation of weak super(sub)-solution that involves the weak derivative of an approximant of $u$.
Let $\rho(x)\in C_{o}^{\infty}(\mathbb{R}^{N})$, $\rho(x)\geqslant 0$ in $\mathbb{R}^{N}$, $\rho(x)\equiv 0$ for
$|x| > 1$ and $\int_{\mathbb{R}^{N}}\rho (x)\,dx=1$, and set
\[
\rho_{h}(x):= h^{-N}\rho\left(x/h\right), \quad
u_{h}(x, t):= h^{-1}\int_{t}\limits^{t+h}\int_{\mathbb{R}^{N}}u(y, \tau)\rho_{h}(x-y)\,dy d\tau.\]

\vskip0.2cm \noindent We fix $t \in (0, T)$ and let $h>0$ be so small that $0<t<t+h<T$.
In \eqref{eq1.4} we take $t_{1}=t$, $t_{2}=t+h$ and replace $\zeta$ by $\int_{\mathbb{R}^{n}}\zeta(y, t)\rho_{h}(x-y)\,dy$.
Dividing by $h$, since the testing function does not depend on $\tau$, we obtain
\begin{equation}\label{eq1.5}
\int_{E\times \{t\}} \left(\frac{\partial u_{h}}{\partial t}\, \zeta+[\mathbb{A}(x, t, \nabla u)]_{h} \nabla \zeta\right)dx
\geqslant\,0\,(\leqslant 0),
\end{equation}
for all $t \in (0\, ,\,  T-h)$ and for all  $\zeta \in W^{1,q}_{o}(E)$, $\zeta \geqslant 0$.  

\subsection{Local Energy Estimates and Critical Mass Lemma}\label{subsec2.1}  \noindent Let $u$ be a weak non-negative super-solution to equation \eqref{eq1.1} in $\Omega_T$, and suppose that for $(\bar{x},\bar{t}) \in \Omega_T$ and $\eta,r>0 $ the following inclusion holds true \[ Q^{+}_{r,\eta}(\bar{x},\bar{t}):=B_{r}(\bar{x})\times(\bar{t}\, ,\,  \bar{t}+\eta)\subset \Omega_T.\] \noindent

\begin{lemma}[Energy Estimates]\label{lem2.1}
 Let $u$ be a non-negative, local weak super-solution to equation \eqref{eq1.1} in $\Omega_T$, and let $\eta,r>0$ and $(\bar{x},\bar{t}) \in \Omega_T$ be as above. For any $\sigma \in (0,1)$, let $\zeta(x,t)=(\zeta_1(x)\zeta_2(t))^q$, for $0 \leq \zeta_i \leq 1$, be a cut-off function such that
 \begin{equation*}
     \begin{cases}
         \zeta_{1} \in C^{\infty}_{o}(B_{r}(\bar{x})): \quad  \zeta_{1}(x)=1 \quad \text{in}\quad  B_{r(1-\sigma)}(\bar{x}), \quad \text{and} \quad \|\nabla  \zeta \|_{\infty}\leq \| \nabla \zeta_{1}\|_{\infty} \leqslant \gamma (\sigma r)^{-1};\\
         \\
                  \zeta_{2} \in C^{1}(\mathbb{R}_{0}^+): \quad \begin{cases} \zeta_{2}(t)=1, \quad t\leqslant \bar{t}+\eta(1-\sigma),\\
         \zeta_{2}(t)=0, \quad  t\geqslant \bar{t}+\eta,\end{cases} \quad  \text{and} \quad \|\partial_t \zeta \|_{\infty}\leq \|\zeta_{2}'\|_{\infty}\leqslant \gamma (\sigma \eta)^{-1}.
     \end{cases}
 \end{equation*}
 \noindent where the $\infty$-norm is taken in $Q_{r, \eta}^+(\bar{x}, \bar{t})$. Let $k$ be any positive constant. Then, if we define  
 \[ [\varphi^{\pm}_{k,r}]=\varphi^{\pm}_{Q^{+}_{r,\eta}(\bar{x},\bar{\eta})}\bigg(\frac{k}{r}\bigg)=\bigg(\frac{k}{r}\bigg)^{p}+ a^{\pm}_{Q^{+}_{r,\eta}(\bar{x},\bar{t})}\bigg(\frac{k}{r}\bigg)^{q},
 \] there exists a positive constant $\gamma$, depending only on the data, such that \begin{multline}\label{eq2.1}
\sup\limits_{\bar{t}<t<\bar{t}+ \eta}\int_{B_{r}(\bar{x})}\zeta(u-k)^{2}_{-} \, dx +
\bigg(\frac{r}{k}\bigg)^{p}\frac{[\varphi^{-}_{k,r}]}{\gamma}\iint\limits_{Q^{+}_{r,\eta}(\bar{x}, \bar{t})}|\nabla [\zeta (u-k)_{-}] |^p\,  dxdt
\\ \leqslant \gamma \sigma^{-q} [\varphi^{+}_{k,r}] \, \bigg(  1+\frac{k^{2}}{\eta [\varphi^{+}_{k,r}]}\bigg) |A^{-}_{k,r,\eta}|,
\end{multline}
\begin{multline}\label{eq2.2}
\sup\limits_{\bar{t}<t<\bar{t}+ \eta}\int_{B_{r}(\bar{x})} \zeta_{1}^{q}(u-k)^{2}_{-}\, dx + \bigg(\frac{r}{k}\bigg)^{p}\frac{[\varphi^{-}_{k,r}]}{\gamma}\iint\limits_{Q^{+}_{r,\eta}(\bar{x}, \bar{t})}|\nabla[\zeta_1^q  (u-k)_{-}]|^{p}\,  dxdt
 \\ \leqslant \int_{B_{r}(\bar{x})\times\{\bar{t}\}} \zeta_{1}^{q} (u-k)^{2}_{-}dx+ \gamma \sigma^{-q} \, [\varphi^{+}_{k,r}]\, |A^{-}_{k,r,\eta}|,
\end{multline}

\end{lemma}

\noindent where $A^-_{k,r,\eta}$ are the $k$ sub-level sets of $u$ in $Q^+_{r,\eta}(\bar{x},\bar{t})$ (see Section \ref{Notation}).\vskip0.1cm \noindent 
Classically, for most parabolic differential equations it is possible to show that the energy estimates, chained with a proper Sobolev-Poincar\'e inequality, imply some sort of measure-theoretical maximum property (see \cite{DiB} for instance). The double-phase equation \eqref{eq1.1} is no exception, provided that a particular geometry is chosen; here we specialize to super-solutions in $Q^+_{4r}$ (see Section \ref{Notation}).

\begin{lemma}[Initial-values Critical Mass] \label{lem2.4}
Let $u$ be a bounded, weak, non-negative super-solution to equation \eqref{eq1.1} in   $Q^{+}_{4r}(\bar{x}, \bar{t})$, with $0\leq u \leq M$. Assume also that for some $0< k < M$
\begin{equation}\label{eq2.7}
u(x,\bar{t})\geqslant k,\quad x\in B_{r}(\bar{x}).
\end{equation}
\noindent Then there exists $\delta \in (0,1)$, depending only on the data, such that for almost all $\displaystyle (x,t)\in  Q^{+}_{r/2,\eta_k}(\bar{x}, \bar{t})$
\begin{equation}\label{eq2.7a}
u(x,t) \geqslant \delta \,k\, ,
\end{equation}
provided that  
\begin{equation}\label{eq2.8}
\eta_k= \frac{k^2}{[\varphi^+_{k,2r}]} \leq (4r)^2 \leq R_o^2 \ , \qquad [\varphi^{+}_{k,2r}]= \bigg( \frac{k}{2r}\bigg)^p + \bigg(\max\limits_{Q_{2r}^+(\bar{x},\bar{t})}a\bigg) \bigg(\frac{k}{2r} \bigg)^q\, \, . 
\end{equation}
\end{lemma}

\begin{proof} 
We consider $(\bar{x},\bar{t})=(0,0)$, to simplify the notation, and an intermediate level $0<\bar{k}<k$. For $n \in \mathbb{N}_0$, we construct 
\[Q_{n}:= Q_{r_n,\eta_k}^+\subset Q_{r, \eta_k}^+=:Q_{o}, \quad \text{being} \quad r_n=r(1+2^{-n})/2,\quad \text{and let} \quad k_n=\bar{k} (1+2^{-n})/2.\] Let us define
\[[\varphi_n^{\pm}]=\bigg( \frac{k_n}{r_n} \bigg)^p + \bigg( a^{\pm}_{Q_{n}} \,    \bigg) \bigg( \frac{k_n}{r_n} \bigg)^q. 
\]\noindent Function $u$ satisfies \eqref{eq2.2} for cut-off functions $\zeta_n=\zeta_1^q$ between $Q_n$ and $Q_{n+1}$ independent of time. The assumption \eqref{eq2.7} simplifies the right-hand side of \eqref{eq2.2} and 
\noindent provides the estimates
\begin{equation}\label{e1}
    \sup_{0<t<\eta_{k}} \int_{B_n} [\zeta_n(u-k_n)_-]^2 \, dx \leq \gamma 2^{nq} [\varphi^+_n] \, \, |[u<k_n]\cap Q_n|,
\end{equation}\noindent and 
\begin{equation}\label{e2}
    \iint_{Q_n} |\nabla [\zeta_n (u-k_n)_-]|^p\, dxdt \leq \gamma 2^{nq} \bigg( \frac{[\varphi^+_{n}]}{[\varphi^-_{n}]} \bigg)\,\bigg(\frac{k_n}{r_n} \bigg)^p\,   |[u<k_n] \cap Q_n|.
\end{equation} \noindent Hence Sobolev's parabolic embedding theorem applies to $[\zeta (u-k_n)_-]$ and \eqref{e1}-\eqref{e2} imply 
\begin{equation*}
    \begin{aligned}
    &(2^{-(n+1)}\bar{k})^p |[u<k_{n+1}] \cap Q_{n+1}| \leq \iint_{Q_{n+1}} (u-k_n)_-^p \, dxdt\\
    & \leq \iint_{Q_{n}} [\zeta_n(u-k_n)_-]^p \, dxdt \\
    & \leq \bigg( \iint_{Q_{n}} [\zeta_n(u-k_n)_-]^{\frac{p(N+2)}{N}}\, dxdt \bigg)^{\frac{N}{N+2}} |[u<k_{n}] \cap Q_{n}|^{\frac{2}{N+2}} \\
    & \leq \gamma \bigg( \sup_{0<t<\eta_k}  \int_{B_{r_n}} [\zeta_n(u-k_n)_-]^2 \, dx \bigg)^{\frac{p}{N+2}} \bigg( \iint_{Q_{n}} |\nabla [\zeta_n (u-k_n)_-]|^p\, dxdt \bigg)^{\frac{N}{N+2}}|[u<k_{n}] \cap Q_{n}|^{\frac{2}{N+2}}\\
    & \leq \gamma 2^{\frac{nq(N+p)}{N+2}}\bigg([\varphi^+_{n}] |[u<k_{n}] \cap Q_{n}| \bigg)^{\frac{p}{N+2}} \bigg( \bigg(\frac{k_n}{r_n} \bigg)^p \frac{[\varphi^+_{n}]}{[\varphi^-_{n}]} \, \, |[u<k_{n}] \cap Q_{n}|\bigg)^{\frac{N}{N+2}}|[u<k_{n}] \cap Q_{n}|^{\frac{2}{N+2}}\\
    & \qquad = \gamma 2^{\frac{nq(N+p)}{N+2}} \bigg( \frac{k_n}{r_n}\bigg)^{\frac{pN}{N+2}}[\varphi^+_{n}]^{\frac{p}{N+2}}  \bigg( \frac{[\varphi^+_{n}]}{[\varphi^-_{n}]}\bigg)^{\frac{N}{N+2}}|[u<k_{n}] \cap Q_{n}|^{1+\frac{p}{N+2}}.
    \end{aligned}
\end{equation*} \noindent
Now we employ condition \eqref{A}, under the assumption $\eta_k \leq  (4 r)^2\leq R_o^2$, therefore we can estimate the ratio $([\varphi^+_n]/[\varphi^-_n])$ with 
\[ [\varphi^+_n] \leq [\varphi^-_n] + A_or_n^{q-p} \bigg( \frac{k_n}{r_n}\bigg)^{q}\leq [\varphi^-_n] \bigg(1+ \frac{A_o k_n^q r_n^{-p}}{(\frac{k_n}{r_n})^p+a^-_{Q_n}(\frac{k_n}{r_n})^q}  \bigg) \leq [\varphi^-_n] \bigg(1+A_oM^{q-p} \bigg), \]
having used also that $k_n \leq k\leq M$. Hence, letting 
\[Y_n= \frac{[u<k_{n}] \cap Q_{n}|}{|Q_n|},\]
and using that $|Q_n|\ge \gamma |Q_{n+1}|$ we obtain 
\begin{equation}
    Y_{n+1} \leq \gamma 2^{\frac{nq(N+p)}{N+2}}  \bigg( \frac{[\varphi^+_{n}] \eta_{k}}{k_n^2} \bigg)^{\frac{p}{N+2}} Y_{n}^{1+\frac{p}{N+2}}\leq \gamma 2^{\frac{nq(N+p)}{N+2}}  \bigg( \frac{2^q [\varphi^+_{k,2r}] \eta_{k}}{(k/2)^2} \bigg)^{\frac{p}{N+2}} Y_{n}^{1+\frac{p}{N+2}} .
\end{equation} \noindent 
We recall here that both $A_o$ and $M$ are structural data. For $0<\delta <1$ to determined, let $\bar{k}= \delta k$. The fast convergence Lemma (see for instance \cite{DiB}, Chap I, Lemma 4.1) gives $Y_{n} \rightarrow 0$ as $n \rightarrow \infty$, provided
\begin{equation}\label{iteration}
Y_0 = \frac{|[u<\delta k]\cap Q_0|}{|Q_0|} \leq  \gamma \bigg( \frac{k^2}{[\varphi^+_{k,2r}] \eta_{k}}\bigg)=: \nu.
\end{equation} \noindent Observe that with our definition of $\eta_{k}$, the number $\nu \in (0,1)$ depends only on the data. In order to prove $Y_0 \leq \nu$, we use again the energy estimates \eqref{eq2.2} to get for $\delta \in (0,1/2)$ the bound
\[
\sup_{0<t<\eta_k} (\delta k)^2 |[u(\cdot, t)<\delta k] \cap B_{2r}|\leq \sup_{0<t<\eta_k} \int_{B_{2r}} (u-2 \delta k)_-^2 \, dx \leq \gamma \, [\varphi^+_{2\delta k,2r}] \, |Q_{2r,\eta}^+| \leq \gamma \delta^p  [\varphi^+_{k,2r}] \, |Q_{2r,\eta}^+|,  \]
where we have used the property $[\varphi^{\pm}_{ck,r}] \leq c^p \, [\varphi^{\pm}_{k,r}]$ for $c \in (0,1)$. Hence
\begin{equation*}
 \begin{aligned}
     Y_0=& \displaystyle{\frac{\int_0^{\eta_{ k}} |[u(\cdot, t)<\delta k] \cap B_{2r}|\, dt}{|Q_0|} } \\
     &\qquad \leq \frac{\eta_{k}\,  \sup_{0<t<\eta_k} |[u(\cdot, t)<\delta k] \cap B_{2r}|}{|Q_0|}\\
     &\qquad \qquad \leq \gamma \frac{\delta^{p-2}}{k^2} [\varphi^+_{k,2r}] \eta_{k} = \gamma \delta^{p-2}/ \nu\quad ,
\end{aligned}
\end{equation*}
\noindent  and condition \eqref{iteration} is satisfied by choosing $\delta$ according to 
\[
Y_0\leq \nu \quad \Leftarrow \quad \delta \leq (\gamma^{-1} \nu^2)^{\frac{1}{p-2}}.
\]
\end{proof}

\begin{remark}
Smaller radii than the levels ensure the previous necessary restriction on $\eta_k$, as
\begin{equation} \label{lungh}
\begin{cases}(k/2r) \ge 1,\\
r\leq R_o/4, \end{cases} 
\qquad \Rightarrow \qquad \begin{cases} \eta_k <(4r)^2,\\
\eta_k \leq R_o^2.
\end{cases}
\end{equation} 
\end{remark}

\vskip0.1cm 
\noindent Now, we need a tool to prolong the information \eqref{eq2.8} to indefinite longer times.\vskip0.1cm 
\noindent Next result roughly states that the estimate \eqref{eq2.7a} is valid for all times that respect the law $|t-\bar{t}| \leq (4r)^2$, at the price of a suitable decay of the level $k$. It is an adaptation of Corollary 3.4 of \cite{GSV} to our double-phase problem.

\begin{corollary}\label{cor} Let the assumptions of Lemma \ref{lem2.4} be satisfied, and suppose the equation \eqref{eq1.1} is satisfied in  $Q^+_{4r}(\bar{x},\bar{t})$, with $0<r<R_o$. Let us define the decreasing function 
\begin{equation*}\label{P}
    \Psi (s)= \frac{s^2}{s^p+a^+_{Q^+_{4r}(\bar{x},\bar{t})} s^q}, \qquad \text{and} \quad \Psi^{-1}\quad \text{its inverse}.
\end{equation*} Then for all $\bar{t}\leq t \leq \bar{t}+ (4r)^2$ and $\delta, \eta_k$ as in \eqref{eq2.8}, the following estimate holds true for all $x \in B_{{r/2}}(\bar{x})$
\begin{equation}\label{DGpower}
u(x,t) \ge \delta k\Psi^{-1} \bigg(1+\frac{(t-\bar{t})}{\eta_k} \bigg). 
\end{equation} 
\end{corollary}

\begin{proof}
\noindent Observe first that, because \eqref{eq2.7} is preserved by diminishing $k$, we can take $0<k<1$. Consider, in the statement of Lemma \ref{lem2.4} the alternatives 
\[  \bar{t} \leq t \leq \bar{t}+ \eta_k \qquad \text{or} \qquad t>\bar{t} +\eta_k. \] \noindent 
In the first case, the application of the aforementioned Lemma turns the information 
\begin{equation*}\label{info}
u(x,\bar{t}) \ge k, \qquad x \in B_{r}(\bar{x}),
\end{equation*} \noindent into 
\[u(x,t)\ge \delta k= \delta k \Psi^{-1} ( \Psi(1)) \ge \delta k \Psi^{-1} (1) \ge \delta k \Psi^{-1}(1 + (t-\bar{t})/\eta_k),\] as both $\Psi$, $\Psi^{-1}$ are decreasing and $\Psi(1)\leq  1$. In the second case, we let 
\[\bar{k}= k \Psi^{-1} \bigg( \frac{t-\bar{t}}{\eta_k}\bigg)  \leq  k,\] and the information 
\[
u(x,\bar{t}) \ge \bar{k} \, \quad \text{in} \quad B_{r}(\bar{x})
\] together again with the use of Lemma \ref{lem2.4} brings us to 
\[u(x,t) \ge \delta \bar{k}, \quad \text{in} \quad B_{r/2}(\bar{x}) \times (\bar{t},\,  \bar{t} + \eta_{\bar{k}}),\] \noindent with 
\[\eta_{\bar{k}} = \Psi\bigg( k\Psi^{-1}\bigg(\frac{t-\bar{t}}{\eta_k}\bigg) \bigg) \ge \Psi(k) \Psi \bigg( \Psi^{-1}\bigg(\frac{t-\bar{t}}{\eta_k}\bigg) \bigg)= \eta_k \bigg(\frac{t-\bar{t}}{\eta_k}\bigg)= (t-\bar{t}).\]
\noindent Here we have used the simple fact that $\Psi(st)\ge \Psi(s) \Psi(t)$ for $s<1$.
\end{proof}

\subsection{Testing with negative powers towards a Reverse H\"{o}lder's inequality}\label{subsec2.5}
\begin{lemma}\label{lem2.2}
Let $(\bar{x},\bar{t}) \in \Omega_T$, and $r, \eta >0$ such that  $Q^{+}_{4r,4\eta}(\bar{x},\bar{t})\subset \Omega_T$. If $u$ is a non-negative, local weak super-solution to equation \eqref{eq1.1} in $\Omega_T$, then for any $\delta\ge 0$, and any $\alpha, \sigma \in (0,1)$, the inequality
\begin{multline}\label{eq2.3}
\frac{1}{1-\alpha}\sup\limits_{\bar{t}<t<\bar{t}+ \eta}\int_{B_{r}(\bar{x})}(u+\delta)^{1-\alpha}\zeta\, dx +\frac{\alpha}{\gamma}\iint\limits_{Q^{+}_{r,\eta}(\bar{x}, \bar{t})}|\nabla [(u+\delta)^{\frac{p-\alpha -1}{p}}{\color{blue}\zeta}]\, |^{p}\,  dxdt + \\+ \frac{\alpha}{\gamma}\iint\limits_{Q^{+}_{r,\eta}(\bar{x}, \bar{t})}a(x,t) |\nabla [(u+\delta)^{\frac{q-\alpha-1}{q}}\zeta]|^q\,   dxdt \leqslant \frac{1}{(1-\alpha)}\|\partial_t \zeta \|_{\infty} \iint\limits_{Q^{+}_{r,\eta}(\bar{x}, \bar{t})}(u+\delta)^{1-\alpha}dxdt +\\ +\gamma \alpha^{1-p} \|\nabla \zeta \|_{\infty}^p \iint\limits_{Q^{+}_{r,\eta}(\bar{x}, \bar{t})}(u+\delta)^{p-\alpha -1}dxdt + \gamma \alpha^{1-q} \|\nabla \zeta \|_{\infty}^q a^{+}_{Q^{+}_{r,\eta}(\bar{x},\bar{t})}\iint\limits_{Q^{+}_{r,\eta}(\bar{x}, \bar{t})}(u+\delta)^{q-\alpha -1} dxdt.
\end{multline} \noindent holds true for any $\zeta_1,\zeta_2$ as in Lemma \ref{lem2.1}, being $\zeta= (\zeta_1\zeta_2)^q$.
\end{lemma}

\noindent The following Lemma constitutes, for nonnegative super-solutions to \eqref{eq1.1}, the reverse H\"older's inequality that we will need for our purpose. 
\begin{lemma}\label{lem2.5}
Let $u$ be a non-negative, bounded, local weak super-solution to equation \eqref{eq1.1} in   $Q^{+}_{r,\eta}(\bar{x}, \bar{t})\subset  Q^+_{r}(\bar{x},\bar{t}) \subset \Omega_T$, with $r<R_o$. Then, for 
all $m\in (0,1)$ and $\delta \ge 0$, there exists a positive constant $\gamma(m)$, depending on the known data and $m$, such that
\begin{equation}\label{eq2.10}
\begin{aligned}
 \frac{1}{r^{p}} \int_{\bar{t}}^{\bar{t}+\eta} \dashint_{B_{r/2}(\bar{x})} &(u+\delta)^{p-2+  \frac{m(p+N)}{N}} \,  dxdt\,\,  + \, \, \frac{  a^+_{Q_{r,\eta}(\bar{x},\bar{t})}}{r^q} \int_{\bar{t}}^{\bar{t}+\eta}  \dashint_{B_{r/2}(\bar{x})} (u+\delta)^{q-2+m(\frac{p+N}{N})} \, dxdt   \\
          &   {}
          \\
          &\leq \gamma(m) I^{m(\frac{p+N}{N})} \bigg{\{} 1+ \eta \bigg( \frac{I^{p-2}}{r^p} + a^+_{Q_{r,\eta}(\bar{x},\bar{t})}\frac{I^{q-2}}{r^q} \bigg) \bigg{\}},
\end{aligned} \end{equation} \noindent where
\begin{equation*}
I:= \sup\limits_{\bar{t}<t< \bar{t}+\eta}\fint\limits_{B_{r}(\bar{x})}u(x,t) dx.
\end{equation*} \noindent The constant $\gamma(m)$ degenerates as soon as $m\downarrow 0$ or $m \uparrow 1$.
\end{lemma} 

\begin{proof}  Let $(\bar{x},\bar{t})$ be the origin (just to ease the notation) and let us define, for $n\in \N \cup \{0\}$,
\[Q_n=  B_n \times (0, \eta), \quad B_n= B_{r_n}, \quad r_n = (r/2)(1+ 2^{-n}), \]
and $\zeta_{n}\in C_o^{1}(B_n)$ a cut-off function such that $\zeta_n \equiv 1$ on $B_{n+1}$, obliged to satisfy
\begin{equation} \label{propzetan} \|\nabla \zeta_n\|_{\infty}:= \|\nabla \zeta_n\|_{L^{\infty}(B_n)}\leq \gamma 2^n/r. \end{equation}
We use H\"older's inequality first with exponent $N/p$ and then with exponent $1/m$ to estimate, in $Q_n$, the quantity
\begin{equation} \label{lasagna}
    \begin{aligned}
\iint_{Q_n}& (u+ \delta)^{p-2+m+\frac{mp}{N}}   \zeta_n^q\, dxdt \\
& \leq \int_0^\eta \bigg( \int_{B_n} (u+\delta)^{m}\, dx \bigg)^{\frac{p}{N}} \bigg( \int_{B_n} [(u+\delta)^{(p-2+m)} \zeta_n^q]^{\frac{N}{N-p}} \, dx \bigg)^{\frac{N-p}{N}}\\
&\leq \int_0^\eta \bigg[ \bigg( \int_{B_n} (u+\delta)\, dx \bigg)^{m} |{B_n}|^{1-m} \bigg]^{\frac{p}{N}} \bigg[ \int_{B_n} \bigg( (u+\delta)^{\frac{p-2+m}{p}} \zeta_n^{\frac{q}{p}} \, dx \bigg)^{\frac{Np}{N-p}} \bigg]^{\frac{N-p}{N}}\\
&\leq \int_0^\eta \bigg[ \bigg( \sup_{0<t<\eta} \dashint_{B_n} u \, dx + \delta \bigg)^m \, |{B_n}|\bigg]^{\frac{p}{N}} \bigg[ \int_{B_n} |\nabla [(u+\delta)^{\frac{p-2+m}{p}}\zeta_n^{\frac{q}{p}}]|^p\, dx\bigg] \, dt,
    \end{aligned}
\end{equation}\noindent by applying Sobolev-Poincar\'e embedding in the last inequality. Now, the first factor of the product on the right-hand side of \eqref{lasagna} is a power of $I$, while we estimate the second integral  on the right-hand side with Lemma \ref{lem2.2} with $m=1-\alpha$ and $\zeta_n=\zeta_1^q$ independent of time, to get from \eqref{lasagna} the inequality
\begin{equation}\label{p}\begin{aligned}
&\iint_{Q_{n+1}} (u+\delta)^{p-2+\frac{m(p+N)}{N}}   \, dxdt \\
 &\leq \iint_{Q_n} (u+\delta)^{p-2+\frac{m(p+N)}{N}}   \zeta_n^q\, dxdt \\
& \leq \gamma(m) (2I^m |{B_n}|)^{\frac{p}{N}} \bigg{\{} \int_{B_n} (u+\delta)^m\, dx+ \iint_{Q_n} \|\zeta_n\|^p_{\infty} (u+\delta)^{p-2+m}+ \|\zeta_n\|_{\infty}^{q} a^+_{Q_n}  (u+\delta)^{q-2+m}\, dxdt \bigg{\}} \\
& \leq \gamma |{B_n}|^{\frac{p}{N}+1} \bigg{\{} I^{m(\frac{p+N}{N})}+ I^{\frac{mp}{N}}\int_0^\eta \dashint_{B_n} \|\zeta_n\|^p_{\infty} (u+\delta)^{p-2+m}+ \|\zeta_n\|_{\infty}^{q} a^+_{Q_n}  (u+\delta)^{q-2+m}\, dxdt \bigg{\}}\\
&\qquad \qquad \qquad \qquad =: \gamma |B_n|^{\frac{p}{N}+1}\, E_n.
\end{aligned}
\end{equation} \noindent  We perform a similar estimate for the phase energy: first we use H\"older's inequality with power $N/p$ and then with $(N-p)/(N-q)$ to get
\begin{equation} \label{b}
    \begin{aligned}
        a^-_{Q_0} \iint_{Q_{n+1}}& (u+\delta)^{q-2+m(\frac{p+N}{N})}\, dxdt \\
        &\leq a^-_{Q_0} \iint_{Q_n} (u+\delta)^{q-2+m(\frac{p+N}{N})} \zeta_n^q\, dxdt \\
         &\leq a^-_{Q_0} \int_0^{\eta} \bigg( \int_{B_n} (u+ \delta)^{m}dx\bigg)^{\frac{p}{N}} \bigg( \int_{B_n} [(u+ \delta)^{q-2+m} \zeta_n^{q}]^{\frac{N}{N-p}} dx \bigg)^{\frac{N-p}{N}}\, dt \\
        & \leq a^-_{Q_0} (2I^m|B_n|)^{\frac{p}{N}} \int_0^\eta \bigg( \int_{B_n} \bigg( [(u+ \delta)^{q-2+m}\zeta_n^q]^{\frac{N}{N-q}} dx\bigg)^{\frac{N-q}{N}} |B_n|^{\frac{q-p}{N}}\, dt\\
        &\leq \gamma  (I^m|B_n|)^{\frac{p}{N}}  |B_n|^{\frac{q-p}{N}}\iint_{Q_n} a(x,t) |\nabla [(u+ \delta)^{\frac{q-2+m}{q}}\zeta_n]|^q\, dxdt\\
        & \leq \gamma (I^m|B_n|^{\frac{q}{N}+1} ) E_n,        
    \end{aligned}
\end{equation} \noindent where we have used again Sobolev-Poincaré inequality in the fourth inequality and Lemma \eqref{lem2.2} in the fifth, denoting the averaged right-hand side of \eqref{eq2.3} in $Q_n$ with $E_n$, as above. Now we use the assumption 
\[ Q_0=Q^+_{r,\eta}\subset Q^+_{r},\]
to apply \eqref{A} and estimate 
\begin{equation*}
    \begin{aligned}
         \frac{a^+_{Q_0}}{r^q}& \int_0^\eta \dashint_{B_n} (u+ \delta)^{q-2+m(\frac{p+N}{N})}\zeta_n^q\, dxdt\\
         &\leq \frac{a^-_{Q_0}}{r^q} \int_0^\eta \dashint_{B_n} (u+ \delta)^{q-2+m(\frac{p+N}{N})}\zeta_n^q\, dxdt+ \frac{AM^{q-p}}{r^p} \int_0^\eta \dashint_{B_n} (u+ \delta)^{p-2+m(\frac{p+N}{N})}\zeta_n^q\, dxdt\\
         &\qquad \qquad \qquad \qquad \leq \gamma I^m (1 +AM^{q-p}) E_n,
    \end{aligned}
\end{equation*} applying \eqref{p}-\eqref{b}. Finally, we estimate $E_n$ by Young's inequality as
\begin{equation}
    \begin{aligned}
    E_n \leq & \, \,   I^{m(\frac{p+N}{N})}+ \int_0^\eta \dashint_{B_n} \|\nabla \zeta_n\|^p_{\infty} (u+ \delta)^{p-2}\bigg( \epsilon (u+ \delta)^{m(\frac{p+N}{N})} + c(\epsilon) I^{\frac{m(p+N)}{N}}\bigg)dxdt + \\
& + \int_0^\eta \dashint_{B_n} \|\nabla \zeta_n\|_{\infty}^{q} a^+_{Q0}  (u+ \delta)^{q-2}\bigg( {\epsilon} (u+ \delta)^{m(\frac{p+N}{N})} + c({\epsilon})  I^{\frac{m(p+N)}{N}}\bigg)\, dxdt\\
&\qquad \leq \gamma \epsilon \int_0^{\eta} \dashint_{B_0} \bigg(\frac{(u+ \delta)^{p-2+m(\frac{p+N}{N})}}{\|\nabla \zeta_n\|^{-p}_{\infty}} + a^+_{Q_0} \frac{(u+ \delta)^{q-2+m(\frac{p+N}{N})}}{\|\nabla \zeta_n\|^{-q}_{\infty}} \,\bigg) dxdt +\\
&\qquad \qquad  \qquad + \gamma I^{m(\frac{p+N}{p})} \bigg(1+ \int_0^{\eta} \dashint_{B_n} \frac{(u+ \delta)^{p-2}}{\|\nabla \zeta_n\|^{-p}_{\infty}}+a^+_{Q_0} \frac{(u+ \delta)^{q-2}}{\|\nabla \zeta_n\|^{-q}_{\infty}}\, dxdt  \bigg) .
    \end{aligned}
\end{equation} \noindent Hence, collecting the terms with $\epsilon$ on a whole initial energetic term $J_0$, and specifying the properties \eqref{propzetan} of $\zeta_n$, we have 
\begin{equation}
    \begin{aligned}
          &J_{n+1}:= \frac{1}{r^p}\int_0^\eta \dashint_{B_{n+1}} (u+ \delta)^{p-2+  \frac{m(p+N)}{N}}  \, dxdt+ \frac{a^+_{Q_0}}{r^{q}\,  } \int_0^\eta \dashint_{B_{n+1}} (u+ \delta)^{q-2+m(\frac{p+N}{N})} \, dxdt \\
          & \leq \gamma I^m E_n
          \\ 
          &\leq \epsilon J_0 + \gamma 2^n I^{m(\frac{p+N}{N})} \bigg{\{} 1+ r^{-p} \int_0^\eta \dashint_{B_n} (u+ \delta)^{p-2}\, dxdt + r^{-q} a^+_{Q_0} \int_0^\eta \dashint_{B_n} (u+ \delta)^{q-2}\, dxdt \bigg{\}}\\
          &\leq \epsilon J_0+\gamma 2^n \bigg{\{} I^{m(\frac{p+N}{N})}+ r^{-p} \int_0^\eta \dashint_{B_n} I^{m(\frac{p+N}{N})} \bigg( (u+ \delta)^{p-2} + r^{-q} a^+_{Q_0}  (u+ \delta)^{q-2}\bigg) dxdt \bigg{\}}\\
           &\leq \epsilon J_0+\gamma 2^n  \bigg{\{} I^{m(\frac{p+N}{N})}+ \epsilon r^{-p} \int_0^\eta \dashint_{B_n} (u+ \delta)^{p-2+m(\frac{p+N}{N})}+ C(\epsilon) I^{p-2+m(\frac{p+N}{N})} \, dxdt +\\
           &\qquad \qquad \qquad \qquad \qquad +r^{-q} a^+_{Q_0} \int_0^\eta \dashint_{B_n} \frac{\tilde{\epsilon}}{2\gamma}  (u+ \delta)^{q-2+m(\frac{p+N}{N})}+ C(\tilde{\epsilon}) I^{q-2+m(\frac{p+N}{N})}\, dxdt \bigg{\}},
    \end{aligned}
\end{equation} \noindent through the use of Young's inequality again, on the last estimate with powers $\frac{p-2+m(\frac{p+N}{N})}{p-2}$ and $\frac{q-2+m(\frac{p+N}{N})}{q-2}$ separately on the terms involving powers of $(u+ \delta)$ and $I$. This finally provides, by choosing again appropriately $\epsilon\in (0,1)$ and reabsorbing the terms in $J_0$, the estimate
\[
J_{n+1}\leq \epsilon J_0+ \gamma \epsilon^{-\gamma} 2^{b n}I^{m(\frac{p+N}{N})} \bigg{\{} 1+ \eta \bigg( \frac{I^{p-2}}{r^p} + a^+_{Q_0} \frac{I^{q-2}}{r^q} \bigg) \bigg{\}}. \]
\noindent Hence a classical iteration provides
\[J_{\infty} \leq \gamma I^{m(\frac{p+N}{N})} \bigg{\{} 1+ \eta \bigg( \frac{I^{p-2}}{r^p} + a^+_{Q_0}\frac{I^{q-2}}{r^q} \bigg) \bigg{\}}.\]
\end{proof}

\subsection{Weak Harnack's Inequality}\label{subsec2.2}
We borrow the following result from \cite{SavSkrYev2}.
\begin{lemma}\label{lem2.3}
Let $u$ be a non-negative, bounded, weak super-solution to equation \eqref{eq1.1} in $Q^{+}_{16r}(\bar{x},\bar{t})$. Then there exist positive numbers $C_{\mathcal{H}}$ and  $b$, depending only on the data, such that 
\begin{equation}\label{eq2.4}
\bar{\mathcal{I}}:= \dashint_{B_{r}(\bar{x})}u(x,\bar{t})dx \leqslant C_{\mathcal{H}} \,\bigg\{r+ r \varphi^{-1}_{Q^{+}_{12r}(\bar{x},\bar{t})}\bigg(\frac{r^{2}}{\eta}\bigg)+ \inf\limits_{B_{2r}(\bar{x})} u(\cdot,t)\,\bigg\},
\end{equation}
for all time levels 
\begin{equation}\label{eq2.5}
\bar{t}+\frac{\eta_{1}}{2}\leqslant t\leqslant \bar{t}+\eta_{1},\quad \eta_{1}:=\min\bigg(\eta\, ,\,  \frac{b r^{2}}{\varphi^+_{Q^{+}_{12r}(\bar{x},\bar{t})}(\frac{\bar{\mathcal{I}}}{r})}\bigg).
\end{equation}
Here $\varphi^{-1}_{Q}(v)$ is the inverse function to the function $\varphi_{Q}^{+}(v):=v^{p-2}+ a^{+}_{Q} \, v^{q-2}$. 
\end{lemma}

\begin{remark}\label{rmk-wh}
\noindent Let $f:\R \rightarrow \R$ be a function that has an increasing inverse $f^{-1}$ and satisfies $ f(\lambda s) \leq \lambda^{q-2} f(s)$ for all $\lambda>1, s \in \R$. By applying $f^{-1}$ to the previous property one gets $\lambda s \leq f^{-1}(\lambda^{q-2}f(s))$ and choosing $s=f^{-1}(x)$ and $\alpha= \lambda^{q-2}$ results in formula 
\[f^{-1}(x) \leq \alpha ^{\frac{-1}{q-2}} f^{-1}(\alpha x), \quad \forall x \in \R, \quad \alpha>1.\]
\noindent The scaling property above translates to $\varphi^{-1}_Q(cx) \leq c^{\frac{1}{q-2}} \varphi^{-1}_Q(x)$ for all $x \in \R$, $0<c<1$.

\end{remark}

\section{Geometric Setting and Auxiliary Results}\label{sect3}
\noindent All the estimates of the previous Sections were of local nature. Here we refine the classic approach to parabolic boundary regularity, in the framework of the double-phase operator \eqref{eq1.1}-\eqref{eq1.2}.
\vskip0.1cm 
\subsection{Preamble}
\noindent Let $(x_o,t_o) \in S_T$ be a point of the lateral boundary of $\Omega_T$. As conditions \eqref{eq1.2} imply $\mathbb{A}(x,t,O)=O$, we extend $\mathbb{A}$ to a vector field $\mathbb{A}(x,t,\xi): \R^{N} \times (0, T] \times \R^N\rightarrow \R^N$ by defining it zero on those vector fields $\xi(x,t): \R^N\times (0, \infty)\rightarrow \R^N$ such that $\xi(x,t)=O$ in the complement of $\Omega_T$. It is easily seen that this extension preserves equation \eqref{eq1.1}-\eqref{eq1.2} in its local definition \eqref{eq1.4}, that now can be formulated in any cylinder 
\[Q_{r}(x_o,t_o)= Q_{r}^-(x_o,t_o) \cup Q_{r}^+(x_o,t_o) \not\subset \Omega_T.\] 
\noindent In this sense we say that some function $v$, that vanishes outside $\Omega_T$, is a local weak sub (super)-solution to \eqref{eq1.1}-\eqref{eq1.2} in such a cylinder. 

\begin{figure}[h]
\centering

\scalebox{0.65}{\begin{tikzpicture}[scale=0.50]

\draw[thick,->] (0,-2) -- (0,25) node[anchor=north west] {\small{$t \in \R$}};
\draw[thick,->] (-2,0) -- (28,0) node[anchor=south west] {\small{$x \in \R^N$}};

\draw[thick] (2,-1) -- (18,-1);

\draw (2,0) rectangle (18,20);

\draw[thick,blue] (11,14) rectangle (25,2);

\draw[thick, red] (7,26) rectangle (27,0.5);

\draw[red] (22,22) node{$Q_{r,r^2}(x_o,t_o)$};

\draw (9,-2) node{$\Omega$};

\draw[thick] (10,14) -- (10,2);

\draw (9,8) node{$\eta$};

\draw[thick] (11,1.5) -- (17.8,1.5);

\draw (16,1) node{$r$};

\draw[blue] (21,5) node{$Q_{r,\eta}^-(x_o,t_o)$};

\draw[blue] (32,5) node{$\eta= \eta^*, \eta_*$};

\draw (14,10) node{$u_k^{\pm}=(u-k^{\pm})_{\pm}$};

\draw (21,10) node{$u_k^{\pm}\equiv 0$};

\draw (18,14) node{$/$};

\draw (20,16) node{$(x_o,t_o)$};

\draw (3,19) node{$\Omega_T$};

\draw (-1,20) node{$T$};

\end{tikzpicture}}
\caption{{\small Scheme of the geometric setting of the proof. For the definition of $\eta^*,\eta_*$ see Subsection \ref{geometric setting} below. Considered a same radius $r$, when $a(x_o,t_o)$ approaches zero, $\eta_*$ stretches to infinity while $\eta^*$ stays unvaried. This motivates the reduction of radii $r<R$ in the former case, according to the size of the phase.}}
 \label{FigA}
\end{figure}
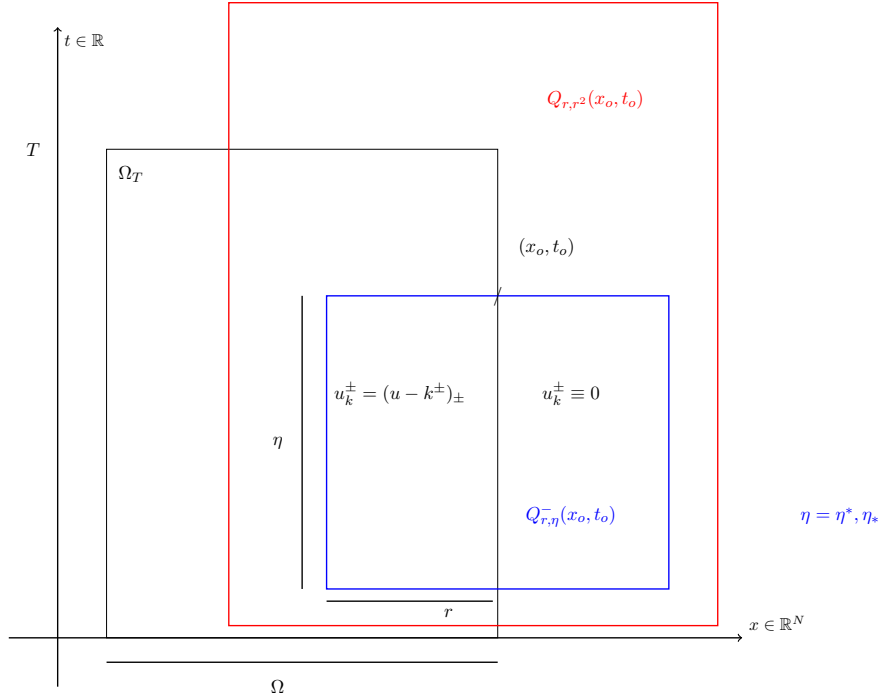

\noindent In the previous Sections we mainly only cared about super-solutions: next Lemma motivates this specialization. Indeed, by extending the equation as above on a cylinder, the truncations $(u-k)_{\pm}$ are sub-solutions, so that $(c-(u-k)_{\pm})$ are non-negative super-solutions, for an appropriate choice of $c>0$. We refer to Lemma 2.1 of \cite{GiaLiaLuk} for more details.

\begin{lemma}\label{lem3.1}
Let $u$ be a local weak solution to equation \eqref{eq1.1}-\eqref{eq1.2} in $\Omega_{T}$ and assume that for a given function $f \in C(\overline{\Omega_T})$ it holds $(u-f)\in V_o^{2,q}$. Let $(x_{o}, t_{o}) \in S_{T}$ and, for some $r>0$ and $0<\eta\leqslant r^{2}$, construct the cylinder \[Q^{-}_{r,\eta}(x_{o},t_{o}):=B_{r}(x_{o})\times(t_{o}-\eta, t_{o})\subset Q_{r}(x_o,t_o).\] We define the zero-extension of the truncations
\begin{equation}\label{trunc1}
   u_{k}^+= \begin{cases} (u-k^+)_{+} & , \ \text{in} \, \, \Omega_T\cap Q_{r,\eta}^-(x_o,t_o)\\
   0 & , \ \text{in} \, \, Q_{r,\eta}^-(x_o,t_o) \setminus \Omega_T 
    \end{cases} \ , \quad \text{for levels} \quad k^+ \ge  \sup_{S_T\cap Q_{r,\eta}^- (x_o,t_o)}f,
\end{equation} \noindent and 
\begin{equation}\label{trunc2}
   u_{k}^-= \begin{cases} (u-k^-)_{-} & , \ \text{in} \, \, \Omega_T\cap Q_{r,\eta}^-(x_o,t_o)\\
   0 & , \  \text{in} \,  \, Q_{r,\eta}^-(x_o,t_o) \setminus \Omega_T   
    \end{cases} \ , \quad \text{for levels} \quad k^- \leq   \inf_{S_T\cap Q^-_{r,\eta}(x_o,t_o)}f.
\end{equation}\noindent Then $u_{k}^{\pm}$ is a weak sub-solution to equation \eqref{eq1.1} 
in $Q^{-}_{r,\eta}(x_{o},t_{o})$. 
\end{lemma} 
\noindent Finally, for $k^{\pm}$ as above, we define $u_{k}^{\pm}$ as in \eqref{trunc1}-\eqref{trunc2}, and we set
\[ w^{\pm}=\mu^{\pm} -u_{k}^{\pm}, \qquad \text{for} \quad \mu^{\pm} \ge \sup_{Q_{r}^{-} (x_o,t_o)} u_k^{\pm}\, \, .\] 
\noindent Evidently $w^{\pm}$ is a  non-negative weak super-solution to equation \eqref{eq1.1} in $Q^{-}_{r}(x_{o},t_{o})$. From here to Section \ref{conclusion}, we drop the superscript $\pm$, because all we need is to work with a generic super-solution $w$.

\subsection{Geometric setting} \label{geometric setting} \noindent The definition of the time-length $\eta>0$ (that must obey $\eta \leq r^2$) needs to distinguish between two different cases: $a(x_{o}, t_{o})=0$ and $a(x_{o}, t_{o}) >0$. 

\subsubsection*{Case $a(x_o,t_o)=0$} \noindent For a number $\gamma^{*}>1$ to be chosen, we let
\[
\delta_{p}(r):=\bigg(\frac{\C_{p}\big(\overline{B_{r}(x_{o})}\setminus\Omega; B_{2r}(x_{o})\big)}{\C_p (\overline{B_r(x_o)};B_{2r}(x_o))}\bigg)^{\frac{1}{p-1}}\]
and consider the time-length
\[\eta^*:= \frac{\gamma^{*}r^{p}}{(\mu \delta_{p}(r))^{p-2}}.\]

\subsubsection*{Case $a(x_o,t_o) >0$}
\noindent Here we set a maximal radius 
\begin{equation} \label{maximal-radius} R^{q-p}:=\dfrac{a(x_{o}, t_{o})}{2A_o},\end{equation} and further we will assume that $r \leqslant \min \{R, R_o\}/24$. This gives us the control on the phase: indeed, in this case 
\begin{equation}\label{eq3.1}
a^{+}_{Q^{-}_{ r,\eta}(x_{o},t_{o})} \leqslant 2 a^{-}_{Q^{-}_{ r,\eta}(x_{o},t_{o})}
\end{equation}
as the simple following computation shows
\begin{equation*}
a^{+}_{Q^{-}_{ r,\eta}(x_{o},t_{o})}-a^{-}_{Q^{-}_{ r,\eta}(x_{o},t_{o})}\leqslant A_o ( r)^{q-p}\leqslant \dfrac{a(x_{o}, t_{o})}{2}\leqslant \dfrac{1}{2}\, a^{+}_{Q^{-}_{ r,\eta}(x_{o},t_{o})}\, \, .
\end{equation*}
\noindent Moreover, the time-length  $\eta_*$ here is defined through the $q$-capacity and the value of $a(x_o,t_o)$, as
\begin{equation*}
\delta_{q}(r):=\bigg(\frac{\C_{q}\big(\overline{B_{r}(x_{o})}\setminus\Omega; B_{2r}(x_{o})\big)}{\C_q(\overline{B_r(x_o)}; B_{2r}(x_o))}\bigg)^{\frac{1}{q-1}},\quad \qquad \eta_*:= \frac{\gamma_{*}r^{q}}{a(x_{o},t_{o})(\mu \delta_{q}(r))^{q-2}}, 
\end{equation*}
\noindent again for a number $\gamma_*>0$ to be chosen (in \eqref{eq3.19}).

\begin{remark} \label{georeq} The conditions $\eta^*, \eta_*<r^2$ imply the estimates
\begin{equation} \label{length-side}\eta^*<r^2 \iff \mu \delta_p(r)\ge (\gamma^*)^{\frac{1}{p-2}}r, \qquad \text{and} \qquad \eta_*<r^2 \iff \mu \delta_q(r)\ge (\gamma_*)^{\frac{1}{q-2}}r.\end{equation}
\end{remark}
\vskip0.3cm \noindent 

\subsection{Capacity estimates} \noindent Now we specialize our estimates towards capacity, considering a test function that vanishes outside a small cylinder. Within the special local geometry chosen, the equation provides a bound for the $p$ or $q$-capacity of $\Omega$ around the point $(x_o,t_o)$, in terms of the averaged $L^1$-norm of $w$. Since a distinction between time lengths is due, for any $0<\eta\leq r^2$ we define 
\begin{equation}\label{I}
    \mathcal{I}(\eta,r)= \sup\limits_{t_{o}-\eta \leqslant t \leqslant t_{o}-\frac{\eta}{4}}\, \,  \fint\limits_{B_{2r}(x_{o})} w(x,t) \, dx, \qquad \text{and} \qquad \mathcal{I}_p(r)= \mathcal{I}(\eta^*,r), \quad \mathcal{I}_q(r)= \mathcal{I}(\eta_*,r).
\end{equation}
\begin{lemma}\label{lem3.2}
Let $u$ be a non-negative, local weak solution of equation \eqref{eq1.1}-\eqref{eq1.2} in $\Omega_T$ and $u=f$ on $S_T$. Fix $(x_o,t_o) \in S_T$ and construct $Q_{r, \eta}^-(x_o,t_o)$ as above (in Section \ref{geometric setting}). There exists a constant $\hat{\gamma}>0$, depending only on the data, such that the following is valid. If $a(x_o,t_o)=0$, then for any $0<r< R_o/2$ we have that
\begin{equation}\label{eq3.2}
\mu\,\delta_{p}(r)\leqslant \hat{\gamma} \mathcal{I}_{p}(r). 
\end{equation}
On the other hand, if $a(x_o,t_o)>0$, for all $0<r<\min\{R_o, R/24\}$, then we find the inequality
\begin{equation}\label{eq3.3}
\mu\,\delta_{q}(r) \leqslant \hat{\gamma} \mathcal{I}_{q}(r) + \hat{\gamma}  \bigg(\frac{r}{R}\bigg)^{q-1}
\frac{1}{(\mu \delta_{q}(r))^{q-2}}.
\end{equation}
\end{lemma}
\begin{proof}

\noindent We divide the argument in two steps: in the first one the special geometry of $\eta^*, \eta_*$ plays no role; while the second one specializes toward $p$ or $q$ capacities. 
\vskip0.2cm  \noindent {\small {\bf STEP 1 - A common potential estimate.}} \vskip0.2cm 

\noindent For any $0<r<\min\{R_o/2, R/16\}$, $0<\eta<r^2$, we construct cylinders $Q_{1}\subset Q_{2} \subset Q_{3}$
\[Q_{1}=B_{r}(x_{o})\times \bigg(t_{o}-\frac{3\eta}{4}, t_{o}-\frac{5\eta}{8}\bigg), \, \,   
Q_{2}=B_{2r}(x_{o})\times \bigg(t_{o}-\frac{7\eta}{8}, t_{o}-\frac{3\eta}{8}\bigg),\, \,  Q_3= B_{4r}(x_o)\times \bigg(t_o-\eta,\, t_o-\frac{\eta}{4}\bigg)\] and let $\zeta \in C^{1}_{o}(Q_{2})$, be a cut-off function between $Q_1$ and $Q_2$, i.e. 
\[\zeta_{|Q_1}\equiv1, \quad \text{and} \quad 0\leqslant \zeta \leqslant 1, \quad |\nabla \zeta| \leqslant \dfrac{2}{r}, \quad 
|\zeta_{t}| \leqslant \dfrac{8}{\eta} \quad \text{in} \,\,  Q_2.\]
\noindent By testing \eqref{eq1.5} with $u_{k,h}\zeta$, for $t\in (t_{o}-\frac{7 \eta}{8}, t_{o}-\frac{3\eta}{8}-h)$, using the fact that $u_{k}$ is a sub-solution of equation \eqref{eq1.1} we obtain
\begin{equation*}
\int_{B_{2r}(x_{o})}\frac{\partial u_{k,h}}{\partial t}u_{k,h}\zeta^{q} dx+
\int_{B_{2r}(x_{o})}[\mathbb{A}(x,t,\nabla u_{k})]_{h} \nabla(u_{k,h}\zeta^{q}) dx \leqslant 0,
\end{equation*}
which yields
\begin{multline*}
\int_{B_{2r}(x_{o})}\frac{\partial w_{h}}{\partial t}w_{h}\zeta dx+
\int_{B_{2r}(x_{o})}[\mathbb{A}(x,t,\nabla u_{k})]_{h} \nabla u_{k,h}\zeta dx \leqslant \\
\leqslant \mu\,\int_{B_{2r}(x_{o})}\frac{\partial w_{h}}{\partial t}\zeta dx+
\int_{B_{2r}(x_{o})}[\mathbb{A}(x,t,\nabla u_{k})]_{h} u_{k,h}\nabla \zeta dx.
\end{multline*}
Now we integrate this inequality over $(t_{o}-\frac{7 \eta}{8}, t_{o}-\frac{3\eta}{8}-h)$. Then, by performing integration by parts in the parabolic terms and finally letting $h \downarrow 0$, while using conditions \eqref{eq1.2}, we find

\begin{equation}
\begin{aligned}
-q \mu \iint_{Q_2} w |\partial_t \zeta | dxdt& - \frac{q}{2} \iint_{Q_2} w^2  |\partial_t \zeta| \, dxdt + \iint_{Q_2} \mathbb{A}(x,t,\nabla u_k) \nabla u_k \zeta  \, dxdt\\
& \leq \mu \iint_{Q_2} \mathbb{A}(x,t,\nabla u_k) \nabla \zeta \, dxdt.
\end{aligned}
\end{equation} \noindent From here, by using the properties of $\zeta$ and the structure conditions \eqref{eq1.2}, we get
\begin{equation}\label{eq3.4}
\iint\limits_{Q_{1}} \varphi(x,t, |\nabla w|) dx dt \leqslant \gamma \mu r^N\mathcal{I}(\eta,r) + \gamma \frac{\mu}{r} \bigg\{\iint\limits_{Q_{2}}|\nabla w|^{p-1} dx dt + \iint\limits_{Q_{2}} a(x,t)|\nabla w|^{q-1} dxdt\bigg\},
\end{equation} for a constant $\gamma>0$ depending only on the data. We abbreviate $\mathcal{I}(\eta,r)=\mathcal{I}$ to ease notation. \newline \noindent Let us estimate the terms on the right-hand side of \eqref{eq3.4}. By H\"older's inequality and being $\mathcal{I}>0$, for any $\bar{m}\in(0,1)$ we obtain
\begin{multline}\label{eq3.5}
\iint\limits_{Q_{2}}|\nabla w|^{p-1} dx dt + \iint\limits_{Q_{2}} a(x,t)|\nabla w|^{q-1} dxdt \leqslant \\ \leqslant
\bigg(\iint\limits_{Q_{2}}(w+\mathcal{I})^{-1-\bar{m}}|\nabla w|^{p} dx dt\bigg)^{\frac{p-1}{p}}\bigg(\iint\limits_{Q_{2}}(w+\mathcal{I})^{(1+\bar{m})(p-1)} dx dt\bigg)^{\frac{1}{p}}+\\+\bigg(\iint\limits_{Q_{2}}a(x,t) (w+\mathcal{I})^{-1-\bar{m}}|\nabla w|^{q} dx dt\bigg)^{\frac{q-1}{q}}\bigg(\iint\limits_{Q_{2}}a(x,t) (w+\mathcal{I})^{(1+\bar{m})(q-1)} dx dt\bigg)^{\frac{1}{q}}.
\end{multline}
Using Lemma \ref{lem2.5} with $m = N(1+\bar{m}(p-1))/(N+p)<1$ we obtain
\begin{equation}\label{eq3.6}
\begin{aligned}
\iint\limits_{Q_{2}}&(w+\mathcal{I})^{(1+\bar{m})(p-1)} dx dt \\
&\leqslant \gamma(\bar{m})r^{N+p}\mathcal{I}^{1+\bar{m}(p-1)}\bigg\{1+\eta\bigg(\frac{\mathcal{I}^{p-2}}{r^{p}}+a^{+}_{Q_{2r}(x_{o},t_{o})}\frac{\mathcal{I}^{q-2}}{r^{q}}\bigg)\bigg\}\\
& 
\\
&=:\gamma(\bar{m})r^{N+p}\mathcal{I}^{1+\bar{m}(p-1)} \, \mathcal{F}(\mathcal{I}).
\end{aligned}
\end{equation}
Similarly, by Lemma \ref{lem2.5} with $m=N(1+\bar{m}(q-1))/(N+p)<1$, we evaluate
\begin{equation}\label{eq3.7} \begin{aligned}
\iint\limits_{Q_{2}}a(x,t) (w+\mathcal{I})^{(1+\bar{m})(q-1)} dx dt \leqslant &a^{+}_{Q_{2r}(x_{o},t_{o})}\iint\limits_{Q_{2}}(w+\mathcal{I})^{(1+\bar{m})(q-1)} dx dt  \\
& \leqslant \gamma(\bar{m})r^{N+q}\mathcal{I}_{1}^{1+\bar{m}(q-1)}\mathcal{F}(\mathcal{I}).
\end{aligned} \end{equation}
Now we use Lemma \ref{lem2.2} for the pair of cylinders $Q_{2}$ and $Q_{3}$, with $\zeta_1 \equiv 1$ on $Q_2$, to compute 
\begin{equation*}
\begin{aligned}
\iint\limits_{Q_{2}}(w+\mathcal{I}_{1})^{-1-\bar{m}}|\nabla w|^{p} dx dt+& \iint\limits_{Q_{2}}a(x,t) (w+\mathcal{I}_{1})^{-1-\bar{m}}|\nabla w|^{q} dx dt  
\\&\leqslant 
\gamma(\bar{m})r^N\mathcal{I}^{1-\bar{m}}  +\frac{\gamma(\bar{m})}{r^{p}}\iint\limits_{Q_{3}}(w+\mathcal{I})^{p-1-\bar{m}} dx dt +\\
&\qquad \qquad \qquad \quad  + a^{+}_{Q_{2r}(x_{o},t_{o})} \frac{\gamma(\bar{m})}{r^{q}}\iint\limits_{Q_{3}}(w+\mathcal{I})^{q-1-\bar{m}} dx dt,
\end{aligned}
\end{equation*}
which by Lemma \ref{lem2.5} with $1-\bar{m}= m(p+N)/N$ yields  the inequality
\begin{equation}\label{eq3.8}
\begin{aligned}
\iint\limits_{Q_{2}}(w+\mathcal{I})^{-1-\bar{m}}|\nabla w|^{p} dx dt+&\iint\limits_{Q_{2}}a(x,t) (w+\mathcal{I})^{-1-\bar{m}}|\nabla w|^{q} dx dt\\
&\leqslant \gamma(\bar{m}) r^{N}\mathcal{I}^{1-\bar{m}}\mathcal{F}(\mathcal{I}).
\end{aligned}
\end{equation}\noindent Collecting  estimates \eqref{eq3.4}--\eqref{eq3.8}, while observing that the powers in \eqref{eq3.5} adjust to 1, we arrive at
\begin{equation}\label{eq3.9}
\iint\limits_{Q_{1}} \varphi(x,t, |\nabla w|) dx dt \leqslant \gamma \mu r^{N}\mathcal{I}\mathcal{F}(\mathcal{I}).
\end{equation}
\noindent Now, as we aim to a capacity estimate, we estimate
\begin{equation} \label{wanted}
\begin{aligned}
    \iint_{Q_2} \varphi(x,t,& |\nabla (\zeta w)|)\, dxdt\\
    &\leq  \iint_{ Q_2} \varphi(x,t, |\nabla w|)\, dxdt +\gamma(C_i) \iint_{ Q_2} \bigg{\{} (w |\nabla \zeta |)^{p}+ a(x,t)  |w\nabla \zeta|^{q} \bigg{\}}\, dxdt,
    \end{aligned}
\end{equation} \noindent where we have used the structure conditions \eqref{eq1.2}. We take care of the second integral term, using Lemma \ref{lem2.5} with $\delta=0$ and $m=N/(p+N)$ to get (as here $\eta<r^2<R_o^2/4$),
\begin{equation} \label{add}
    \begin{aligned}
        \iint_{Q_2} w^p|\nabla \zeta|^p\, dxdt +\iint_{Q_2} a(x,t) & w^q|\nabla \zeta|^q\, dxdt\\
        & \leq \frac{\mu}{r^p} \iint_{Q_3} w^{p-1}\, dxdt+ a^+_{Q_2}\frac{\mu}{r^q} \iint_{Q_3} w^{q-1} \, dxdt\\
        & \qquad \qquad \leq \gamma \mu r^N \mathcal{I} \mathcal{F}(\mathcal{I}).
    \end{aligned}
\end{equation} \noindent Hence finally, joining estimates \eqref{eq3.9} and \eqref{add} into \eqref{wanted} we obtain the potential estimate
\begin{equation}
    \label{potentissima}
    \iint_{ Q_2} \varphi(x,t, |\nabla (\zeta w)|) \, dxdt \leq \gamma \mu r^N \mathcal{I}(\eta,r) \mathcal{F}(\mathcal{I}(\eta,r)).
\end{equation}

\vskip0.1cm  \noindent {\small {\bf STEP 2 - Geometry enters into play.}} \vskip0.1cm 
\noindent Here we divide the study in two cases depending on the value of the phase at the boundary point. 

\noindent If $a(x_{o},t_{o})=0$, we fix $\eta=\eta^*$ as above (Section \ref{geometric setting}) and proceed by contradiction: we assume that for any $\varepsilon\in(0,1)$ (to be determined and depending only on the data) the converse inequality 
\begin{equation}\label{eq3.10}
\mathcal{I}_{p}\leqslant \varepsilon \mu \delta_{p}(r)
\end{equation}
holds true, because otherwise inequality \eqref{eq3.2} is found. Now, by the definition, the scaling properties of the $p$-capacity $\C_{p}(B_{r}(x_{o}); B_{2r}(x_{o}))= \gamma r^{N-p}$ for a positive constant $\gamma$ depending only on $N$ and $p$, and our choice of $\eta^*$, we have 
\begin{equation}\label{eq3.11}
\iint\limits_{Q_{2}} \varphi(x,t, |\nabla  (\zeta w)|) dx dt\geqslant \frac{3}{4} \mu^{p}\eta^* \C_{p}(B_{r}(x_{o})\setminus \Omega; B_{2r}(x_{o}))\ge 
 \gamma \gamma^{*} \mu^{2} \delta_{p}(r) r^N .
\end{equation}
Moreover, since $0<\eta<r^2<R_o^2$ condition \eqref{A} is in force and
\begin{equation*}
a^{+}_{Q_{2r}(x_{o},t_{o})}\dfrac{\mathcal{I}_p^{q-2}}{r^{q}}\leqslant A_o (2r)^{q-p}\dfrac{\mathcal{I}_p^{q-2}}{r^{q}}\leqslant \gamma \dfrac{\mathcal{I}_p^{p-2}}{r^{p}}, \quad \text{with} \quad \gamma=A_o (2M)^{q-p},
\end{equation*}\noindent so that the inequalities  \eqref{potentissima}-\eqref{eq3.11}, chained, can be rewritten as
\begin{equation}\label{eq3.12}
\mu \delta_{p}(r) r^{N} \leqslant \frac{\gamma}{\gamma^{*}} r^{N}\mathcal{I}_p +\frac{\gamma}{\gamma^{*}} r^{N-p}\eta^* \mathcal{I}_p^{p-1}\leqslant \gamma(\varepsilon +\varepsilon^{p-1}) \mu \delta_{p}(r) r^{N}.
\end{equation}
Choosing $\varepsilon$ small enough, such that $\gamma(\varepsilon +\varepsilon^{p-1})=\frac{1}{2}$, a contradiction to  \eqref{eq3.10} is reached. This proves inequality \eqref{eq3.2}.\vskip0.2cm 

\noindent Now if $a(x_{o},t_{o}) >0$, we fix $\eta=\eta_*$ for $0<r<R/16$ (still referring to Section \ref{geometric setting}) and we assume that for any  $\epsilon \in (0,1)$, the estimate 
\begin{equation}\label{eq3.13}
\mathcal{I}_{q} + \bigg(\frac{r}{R}\bigg)^{q-1}
\frac{1}{(\mu \delta_{q}(r))^{q-2}} \leqslant \varepsilon \mu \delta_{q}(r)
\end{equation}
holds true; otherwise, inequality \eqref{eq3.3} would be in force. By the assumption $0<r<R/16$ the estimate \eqref{eq3.1} is valid, while the definition of $q$-capacity and the choice of $\eta_*$ imply
\begin{equation}\label{eq3.14}
\iint\limits_{Q_{2}} \varphi(x,t, |\nabla (\zeta w)|) dx dt\geqslant  a(x_{o}, t_{o}) \mu^{q}\  \frac{3}{4} \eta_* \ \C_{q}(B_{r}(x_{o})\setminus \Omega; B_{2r}(x_{o})) \ge   \gamma \gamma_{*} \,  \mu^{2} \delta_{q}(r) r^{N},
\end{equation}
and \eqref{potentissima}-\eqref{eq3.14}, chained, can be rewritten as
\begin{equation}\label{eq3.15}
\begin{aligned} 
\mu \delta_{q}(r) r^{N}& \leqslant \frac{\gamma}{\gamma_{*}} r^{N} \mathcal{I}_q + \frac{\gamma}{\gamma_{*}} a(x_{o}, t_{o}) \eta_* r^{N-q} \mathcal{I}_q^{q-1}+ \frac{\gamma}{\gamma_{*}} \eta_* r^{N-p} \mathcal{I}_q^{p-1}\\
& \leqslant \gamma( \varepsilon +\varepsilon^{q-1}) \mu \delta_{q}(r) r^{N}+ \frac{\gamma}{\gamma_{*}} \varepsilon^{p-1} \mu^{p-1} \eta_* \delta_{q}(r)^{p-1} r^{N-p}\\
&\leqslant \gamma( \varepsilon +\varepsilon^{p-1}+ \varepsilon^{q-1}) \mu \delta_{q}(r) r^N+ \gamma \varepsilon^{p-1} \frac{r^{N+q-1}}{a(x_{o},t_{o})^{\frac{q-1}{q-p}}(\mu \delta_{q}(r))^{q-2}}\\
&=\gamma( \varepsilon +\varepsilon^{p-1}+ \varepsilon^{q-1}) \mu \delta_{q}(r) r^{N}+\gamma \varepsilon^{p-1}\bigg(\frac{r}{R}\bigg)^{q-1}
\frac{r^{N}}{(\mu \delta_{q}(r))^{q-2}},
\end{aligned}
\end{equation} \noindent where in the second inequality we have used \eqref{eq3.13} as $\mathcal{I} \leq \varepsilon \mu \delta_q(r)$ and the definition of $\eta_*$, while in the third inequality we have used Young's inequality with $(q-1)/(p-1)$ and its conjugate $(q-1)/(q-p)$ weighted on $\mu \delta_q(r) r^N$, separating the term $\varepsilon^{p-1}$ from the remainder. To arrive to the wanted contradiction, it is enough to choose $\varepsilon$ such that $\gamma( \varepsilon +2\varepsilon^{p-1}+ \varepsilon^{q-1})=1/2$. This completes the proof of Lemma \ref{lem3.2}.
\end{proof}

\begin{lemma}\label{lem3.3} Let the assumptions of Lemma \ref{lem3.2} be valid. Then, in the case $a(x_o,t_o)=0$, there exists a constant $C_{p}>0$, depending only on the data, such that, either
\begin{equation}\label{alternativep}
    \mu \delta_p(r) \leq 2C_p r,
\end{equation}\noindent or 
\begin{equation}\label{eq3.22}
\sup\limits_{Q^{-}_{\frac{r}{2},\frac{\eta}{8}}(x_{o},t_{o})}u_{k}\leqslant \mu\,\bigg(1-\frac{1}{2C_{p}}\delta_{p}(r)\bigg),
\end{equation} \noindent In the case $a(x_o,t_o)>0$, there exists a constant $C_q>0$, depending only on the data, such that either 
\begin{equation}\label{alternativeq}
    \mu \delta_q(r)\leq 4 C_q  r+ ( 4C_{q})^{\frac{1}{q-1}} \bigg(\frac{r}{R} \bigg),
\end{equation} or
\begin{equation}\label{eq3.23}
\sup\limits_{Q^{-}_{\frac{r}{2},\frac{\eta}{8}}(x_{o},t_{o})}u_{k}\leqslant \mu\,\bigg(1-\frac{1}{2C_{q}}\delta_{q}(r)\bigg).
\end{equation}
\end{lemma}

\begin{proof}
Referring Section \ref{geometric setting} and Lemma \ref{lem3.2}, we let $\eta =\eta^*/2, \eta_*/2$ in the two cases, and considering the continuity of the function 
\[
[t_o-\eta, \, t_o] \,  \ni t \, \rightarrow \, \fint_{B_{2r}(x_o)} w\, dx,
\]
we let $t_{1} \in [t_{o}-\eta, t_{o}-\eta/2]$ be the point such that inequality \eqref{eq3.2} (or \eqref{eq3.3})) is achieved, i.e.
\begin{equation}\label{t1}  \mathcal{I}_{1}=\sup\limits_{t_{o}-\eta \leqslant t \leqslant t_{o}-\eta/2} \, \, \fint\limits_{B_{2r}(x_{o})} w \, dx=\fint\limits_{B_{2r}(x_{o})} w(x,t_{1})\,  dx,\end{equation} \noindent depending on the choice of $\eta$.  Now we apply the weak Harnack inequality \ref{eq2.4} with \[\eta_{1}=\min \bigg(\dfrac{\eta}{4}\, ,\, \,  \dfrac{4b r^{2}}{\varphi^+_{Q_{ 24r}^+(x_{o},t_{1})}(\frac{\mathcal{I}_{1}}{2r})}\bigg),\] 
which yields
\begin{equation}\label{eq3.18}
\mathcal{I}_{1} \leqslant \gamma \bigg\{r + r \varphi^{-1}_{Q^{+}_{24 r}(x_{o},t_{1})}\bigg(\frac{2r^{2}}{\eta}\bigg)+ \inf\limits_{B_{ 4r}(x_{o})}w(\cdot, t_{2})\bigg\}, \qquad \text{at} \quad  t_{2} = t_{1}+\eta_1/2<t_o-\eta/4.
\end{equation} \noindent We want to estimate the second term on the right-hand side of \eqref{eq3.18}: to this aim, we apply Remark \ref{rmk-wh}. If $a(x_{o},t_{o})=0$, then we evaluate
\begin{multline*}
r \varphi^{-1}_{Q^{+}_{24 r}(x_{o},t_{1})}\bigg(\frac{2 r^{2}}{\eta^*}\bigg) \leqslant  {\gamma}(\gamma^{*})^{\frac{-1}{q-2}} r\varphi^{-1}_{Q^{+}_{24 r}(x_{o},t_{1})}\bigg(\bigg(\frac{\mu\delta_{p}(r)}{r}\bigg)^{p-2}\bigg)\leqslant \\ \leqslant \gamma (\gamma^{*})^{\frac{-1}{q-2}} r \varphi^{-1}_{Q^{+}_{24 r}(x_{o},t_{1})}\bigg(\varphi_{Q^{+}_{24 r}(x_{o},t_{1})}\bigg(\frac{\mu\delta_{p}(r)}{r}\bigg)\bigg)= {\gamma}
(\gamma^{*})^{\frac{-1}{q-2}} \mu \delta_{p}(r).
\end{multline*}
Similarly, if $a(x_{o}, t_{o})>0$, we use the condition $\eta_*<r^2$ to get $|t_o-t_1|\leq \eta_* < r^2/2$ and similarly to \eqref{eq3.1} we can estimate 
\[ a^+_{Q^+_{r}(x_o,t_1)} \leq 2 a^{-}_{Q_{r}(x_o,t_1)} \leq 2 a(x_o,t_o) \leq 2 a^+_{Q^+_{r}(x_o,t_1)}, 
\] because 
\[(x_o,t_o) \in Q^+_{r}(x_o,t_1) \subseteq Q_{2r}(x_o,t_o).\]
\noindent Hence using this fact in the second inequality, a similar computation yields
\begin{multline*}
r \varphi^{-1}_{Q^{+}_{{24 r}(x_{o},t_{1})}}\bigg(\frac{2 r^{2}}{\eta_*}\bigg) \leqslant {\gamma}\gamma_{*}^{-\frac{1}{q-2}} r\varphi^{-1}_{Q^{+}_{24 r}(x_{o},t_{1})}\bigg(a(x_{o},t_{o})\bigg(\frac{\mu\delta_{q}(r)}{r}\bigg)^{q-2}\bigg)\leqslant \\ \leqslant {\gamma}\gamma_{*}^{-\frac{1}{q-2}} r \varphi^{-1}_{Q^{+}_{24 r}(x_{o},t_{1})}\bigg(\varphi_{Q^{+}_{24 r} (x_{o},t_{1})}\bigg(\frac{\mu\delta_{q}(r)}{r}\bigg)\bigg)\leqslant {\gamma} \gamma_{*}^{-\frac{1}{q-2}} \mu \delta_{q}(r).
\end{multline*}
From this and \eqref{eq3.18}, using Lemma \ref{lem3.2} and choosing $\gamma_{*}, \gamma^*$ by the conditions 
\begin{equation}\label{eq3.19}
\gamma_{*}= \gamma^{*}=(2\gamma )^{q-2},
\end{equation}
we arrive at 
\begin{equation}\label{eq3.20}
\mu \delta_{p}(r)\leqslant \gamma \big\{r + \inf\limits_{B_{4r}(x_{o})}w(\cdot, t_{2})\big\},\quad\text{if}\quad a(x_{o},t_{o})=0,
\end{equation}
and
\begin{equation}\label{eq3.21}
\mu\delta_{q}(r)\leqslant \gamma \big\{r+\bigg(\frac{r}{R}\bigg)^{q-1}\frac{1}{(\mu \delta_{q}(r))^{q-2}}+\inf\limits_{B_{4 r}(x_{o})}w(\cdot, t_{2})\big\},\quad \text{if}\quad a(x_{o}, t_{o})>0,
\end{equation}
for $t_{2}=t_1+\eta_1/2<t_o-\eta/4$ by our choice of $\eta_1$. \vskip0.1cm \noindent  We observe that at the moment the quantitative location of $t_2$ is undetermined, because of the unknown $t_1$ (see Figure \ref{FIG2}). To complete the proof, we use Corollary \ref{cor} for the function $w$, as
\begin{equation*}
    \begin{cases} 
w(x,t_2)\ge k_p= {\gamma}^{-1} \mu\delta_p(r) -r,   & x \in  B_{2r}(x_o), \quad \text{if} \quad a(x_o,t_o)=0, \\
w(x,t_2)\ge k_q={\gamma}^{-1} 
 \mu\delta_q(r) - r - ({r}/{R})^{q-1} \frac{1}{(\mu \delta_q (r))^{q-2}},  &x\in B_{2r}(x_o), \quad  \text{if} \quad a(x_o,t_o)>0,
\end{cases}
\end{equation*}
where $k_p,k_q$ are positive by assumptions \eqref{alternativep}-\eqref{alternativeq}. As we have 
\[t_1<t_o-\eta/2, \quad   \text{for} \quad \eta=\eta^*, \eta_*,\] then Corollary \ref{cor} implies that 
\begin{equation*}
    \begin{cases} 
w(x,t)\ge \delta^*(t) k_p, \qquad \quad  (x,t) \in  B_{{\color{blue}r}}(x_o) \times (t_2, t_2+\eta_{k_p}), & \text{if} \quad a(x_o,t_o)=0, \\
w(x,t)\ge \delta_*(t) k_q, \quad \qquad  (x,t) \in  B_{r}(x_o) \times (t_2, t_2+\eta_{k_q}), & \text{if} \quad a(x_o,t_o)>0,
\end{cases}
\end{equation*} for $\eta_{k_p}, \eta_{k_q}$ referred to levels $k_p,k_q$ as given in \eqref{eq2.8}, and
\[\delta^*(t)= \delta \Psi^{-1}\bigg(1+\frac{t-t_{2}}{\eta_{k_{p}}} \bigg), \qquad \quad \delta_*(t)= \delta \Psi^{-1}\bigg(1+\frac{t-t_{2}}{\eta_{k_{q}}} \bigg).\]
\noindent So we move the point-wise information on $w$ from $t_1$ to $t_o$, hence travelling a distance smaller than $\eta$:
\[\delta^*(t)\ge  \delta \Psi^{-1}\bigg(1+\frac{\eta}{\eta_{k_{p}}} \bigg)=\Psi^{-1}\bigg(1+\frac{\gamma^*[\mu(\delta_p(r)]^{2-p} r^p}{\Psi(\hat{\gamma}^{-1} \mu \delta_p(r)-r)} \bigg) \ge \Psi^{-1}\bigg(1+\frac{\Psi(\gamma^*\mu(\delta_p(r))}{\Psi(\hat{\gamma}^{-1} \mu \delta_p(r))} \bigg) =: C_p^*,\]

\[\delta_*(t)\ge  \delta \Psi^{-1}\bigg(1+\frac{\eta}{\eta_{k_{q}}} \bigg)=\Psi^{-1}\bigg(1+\frac{\gamma_*[\mu(\delta_q(r)]^{2-q} r^q}{\Psi(\hat{\gamma}^{-1} \mu \delta_q(r)-r)} \bigg) \ge \Psi^{-1}\bigg(1+\frac{\Psi(\gamma_*\mu(\delta_q(r))}{\Psi(\hat{\gamma}^{-1} \mu \delta_q(r))} \bigg) =: C_q^*.\]
\noindent This implies, in the case $a(x_o,t_o)=0$ and \eqref{alternativep} violated, the estimate
\begin{equation}\label{eq3.16}
\mu \delta_{p}(r) \leqslant C_{p}^* \bigg(\mu-\sup\limits_{Q^{-}_{\frac{r}{2},\frac{\eta}{8}}(x_{o},t_{o})}u_{k}\bigg)+ C_{p}^* r.
\end{equation} Similarly, in the case $a(x_o,t_o)>0$ and \eqref{alternativeq} violated, the above procedure ensures
\begin{equation}\label{eq3.17}
\mu \delta_{q}(r) \leqslant C_{q}^* \bigg(\mu-\sup\limits_{Q^{-}_{\frac{r}{2},\frac{\eta}{8}}(x_{o},t_{o})}u_{k}\bigg)+C_{q}^* r+ C_{q} \bigg(\frac{r}{R}\bigg)^{q-1}\frac{1}{(\mu \delta_{q}(r))^{q-2}}.
\end{equation} \noindent The conclusion follows therefore by implementing the assumption that \eqref{alternativep}-\eqref{alternativeq} are violated into the estimates \eqref{eq3.16}-\eqref{eq3.17} above. 

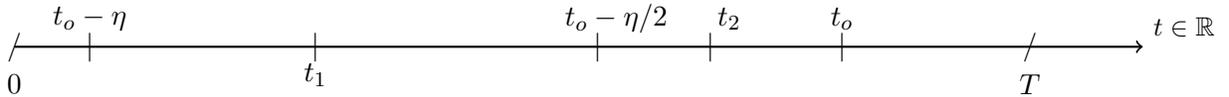
\begin{figure}[h] \label{FIG2}
\centering

{\begin{tikzpicture}[scale=0.50]

\draw[thick,->] (-2,0) -- (28,0) node[anchor=south west] {\small{$t \in \R$}};

\draw (0,0.75) node{$t_o-\eta$};

\draw (6,-0.75) node{$t_1$};


\draw (14,0.75) node{$t_o-\eta/2$};

\draw (17,0.75) node{$t_2$};

\draw (20,0.75) node{$t_o$};

\draw (-2,-1) node{$0$};

\draw (25,-1) node{$T$};

\draw (-2,0) node{$/$};

\draw (25,0) node{$/$};

\draw (0,0) node{$|$};

\draw (6,0) node{$|$};

 \draw (16.5,0) node{$|$};

\draw (13.5,0) node{$|$};

\draw (20,0) node{$|$};

\end{tikzpicture}
\caption{{\small Comparing time lengths in proof of Lemma \ref{lem3.3}}.}
 \label{FigA}}
\end{figure}
    
\end{proof}

\section{Proof of Theorem \ref{boundary-estimate}}\label{conclusion}
 \noindent We begin with a preliminary consideration. The divergence of Wiener's integral at $(x_o,t_o)$ implies that there exists a suitable sequence of radii that allows to apply Lemma \ref{lem3.3} iteratively. 
\begin{lemma}\label{itera} Let $p>1$, $\mu_o>0$ and $\bar{C}, C_1 >1$ be given numbers. Assume that for a certain $\rho_o>0$ it holds both
\begin{equation}\label{W}
    \int_0^{\rho_o} \delta_p (\rho) \frac{d\rho}{\rho}= \infty,
\end{equation} and 
\[\mu_o \delta_p(\rho_o) \ge \bar{C} \rho_o.\]
\noindent Then, for any $\tilde{\gamma}>0$ fixed, there exists a decreasing sequence of radii $\{ \rho_j\}_{j \in \N}$ such that, by defining
\[\mu_j= (1-1/(2C_1))\mu_{j-1}, \qquad \mu_0=\mu_o,\quad \rho_0=\rho_o\, ,\] it has the following properties for all $j \in \N\cup \{0\}$ :
\begin{equation}\label{1}
    \mu_j \delta_p(\rho_j) \ge \bar{C} \rho_j;
\end{equation}

\begin{equation}\label{2}
  2 \eta_{j+1}:= 2 \tilde{\gamma} \rho_{j+1}^p (\mu_{j+1} \delta_p(\rho_{j+1}))^{2-p} \leq  \tilde{\gamma}\rho_{j}^p (\mu_{j} \delta_p(\rho_{j}))^{2-p}= \eta_{j};
\end{equation}

\begin{equation}\label{3}
\forall l \in \N\, \,  \exists n(l)\in \N : \qquad  
    \sum_{j=0}^{l-1} \delta_p(\rho_j) \ge \frac{1}{\gamma_3} \sum_{i=0}^{n(l)} \delta_p(\sigma^{i}\rho_0) \ge \frac{1}{\gamma_4}\int_{\rho_l}^{\rho_0} \frac{\delta_p(\rho) \, d\rho}{\rho}\, .
\end{equation}
    
\end{lemma}

\begin{remark}\label{rmkq}
    
    \noindent If in the previous Lemma we choose $\bar{C}= \tilde{C} (1+R^{-1}+ R^{\frac{p-q}{q-2}})$, with the choice $\rho_o\leq R=({a(x_o,t_o)}/{2}A_o)^{\frac{1}{q-p}}$, then \eqref{1}-\eqref{2}-\eqref{3}
 hold true for the exponent $q$ instead of $p$ and condition \eqref{1} is replaced by 
 \begin{equation}\label{raggi} \mu_j\delta_q(\rho_j) \ge \tilde{C} (1+R^{-1}+R^{\frac{p-q}{q-2}}) \rho_j\, \, .\end{equation}
 Regarding the three terms on the right-hand side of \eqref{raggi}: the first and second one are linked to the requirement of \eqref{alternativeq}; while the third one is given to free the choice of $\tilde{C}$ from $a(x_o,t_o)$ when requiring in Remark \ref{georeq},
 \[Q_{\rho_o, \eta^*} \subseteq Q_{\rho_o, \rho_o^2}, \quad \text{with} \quad  \eta^*= \frac{\gamma_* \rho_0^q}{a(x_o,t_o) (\mu(\rho_0)\delta_q(\rho_o))^{q-2}}\, \, .\]
\end{remark} 

\noindent Lemma \ref{itera} is an adaptation to our framework of an argument of the capacitary estimate between integral and sum of \cite{GiaLiaLuk}, while the extraction of the sequence is modeled after \cite{Skr1}. The novelty is that we assume a priori that $\rho_o$ satisfies \eqref{1}, and we extract the sequence $\{\rho_j\}_{j \in \N}$ starting from $\rho_o$.

\begin{proof} Let $\eta_0= \tilde{\gamma} \rho_0^p(\mu_o\delta_p(\rho_0))^{2-p}$, and 
if 
\[\delta_p(\rho_0(1-1/(2C_1))/2)\ge 1/2 \ \delta_p(\rho_0), \quad \Rightarrow \quad \text{and set} \quad \rho_1= \sigma \rho_0, \quad \quad \sigma:= (1-1/(2C_1))/2;\] \noindent while if otherwise 
\[\delta_p(\rho_0(1-1/(2C_1))/2)< 1/2 \  \delta_p(\rho_0),\] let $i_1\in \N$ be the smallest number such that 
\begin{equation} \label{4} \delta_p(\sigma^{i_1} \rho_0) \ge 2^{-i_1}\delta_p(\rho_0), \quad \Rightarrow \quad \text{and set} \quad  \rho_1= \sigma^{i_1} \rho_0.\end{equation} The choice of $i_1$ is possible, being otherwise 
\[\delta_p(\rho_0(1-1/(2C_1))/2)< 2^{i} \delta_p(\rho_0) \quad \forall i \in \N \qquad \Rightarrow \quad \sum_{i\in \N} \frac{\delta_p(\sigma^{i}\rho_0)}{2^i} < \infty,\] that is a contradiction with \eqref{W}, because of the property
\begin{equation*}
    \begin{aligned}
        \bigg[ \frac{\mathcal{C}_p(B(x_o, 2^{-(k+1)}\rho_0); B(x_o,\rho_0))}{\gamma 2^{-k(N-p)}\rho_0}\bigg]^{\frac{1}{p-1}} \ln (2) &\leq \int_{2^{-(k+1)}\rho_0}^{2^{-k}\rho_0}  \bigg[ \frac{\mathcal{C}_p(B(x_o,t); B(x_o,\rho_0))}{t^{(N-p)}}\bigg]^{\frac{1}{p-1}} \frac{dt}{t}\\
        & \leq    \bigg[ \frac{\mathcal{C}_p(B(x_o, 2^{-k}\rho_0); B(x_o,\rho_0))}{\gamma 2^{-(k+1)(N-p)}\rho_0}\bigg]^{\frac{1}{p-1}} \ln (2).
    \end{aligned}\end{equation*}

\noindent Henceforth, we have \eqref{1} for $i=1$, which is 
\[\mu_1 \delta_p(\rho_1) =\mu_o(1-1/(2C_1))\delta_p(\rho_1) \ge \frac{\mu_o}{2^{i_1}}\delta_p(\rho_0)\ge  \frac{\bar{C} \rho_0}{2^{i_1}}\ge \bar{C} \rho_1.\]
\noindent Point \eqref{2} follows from \eqref{4} and the definition of $\eta_1, \eta_0$, with a simple computation. Finally, the choice of $i_1$ to be the smallest number is finally useful to have \eqref{3}, as  
\[ \sum_{i=0}^{i_1}\ \delta_p(\sigma^i \rho_0) \leq \delta_p(\rho_0) + 2 \delta_p (\rho_0) +\delta_p(\sigma^{i_1}\rho_0) \leq 3 \delta_p(\rho_0) + \delta_p(\rho_1) \leq 3 \sum_{i=1,2} \delta_p(\rho_i), \]
where (in the first inequality) we have contradicted condition \eqref{4} for all previous $0<i<i_1$ to get 
\[2 \delta_p(\rho_0) \ge \sum_{i=1}^{i_1-1} \delta_p(\sigma^i \rho_0).\]
By induction, we suppose the statement of the Lemma to be valid until step $n$ and we prove it for $(n+1)$. Thus, we choose $\rho_{n+1}=\sigma^{i_{n+1}-i_n}\rho_0$ for $i_{n+1}$ being the smallest number in $\{i_{n}, i_{n}+1, \dots\}$ satisfying as before \eqref{4} with $i_{n+1}$ instead. Then all the argument flows in the same style until we arrive to condition \eqref{3}: here we contradict assumption \eqref{4} for all previous $i \in \{i_n, \dots, i_{n+1}\}$, and use the inductive hypothesis to obtain 
\begin{equation}
    \begin{aligned}
        \sum_{i=0}^{i_{n+1}} \delta_p(\sigma^i\rho_0) &\leq  \sum_{i=0}^{i_{n-1}} \delta_p(\sigma^i\rho_0) + \delta_p(\rho_n)+\sum_{i=i_n+1}^{i_{n+1}-1} \delta_p(\sigma^i\rho_0) + \delta_p(\rho_{n+1})\\
        &\leq 3 \sum_{j=0}^{n-1} \delta_p(\rho_j) + \delta_p(\rho_n) + \delta_p(\rho_n) \sum_{i=i_n+1}^{i_{n+1}-1} 2^{i_n-i} + \delta_p(\rho_{n+1})\\
        &\leq 3 \sum_{j=0}^{n+1} \delta_p(\rho_j).
    \end{aligned}
\end{equation} \noindent This ensures that the first inequality of \eqref{3} occurs, with $\gamma_3=3$ and $n(l)=i_l$. \newline
\noindent To prove the third one, we follow a reasoning similar to \cite{GiaLiaLuk}. For a condenser $(K,B_{2r})$, Lemma 2.16 of \cite{HeiKipMar} states that for $p>1$ and when $0<r\leq s\leq 2r$, then there exits $\gamma(s,N)>0$ such that 
\[\frac{1}{\gamma} \C_p(K; B_{2r}) \leq \C_p(K, B_{2s}) \leq \gamma \C_p(K; B_{2r}).\]\noindent Hence, by the previous consideration and the monotonicity of the capacity in the first argument, we obtain \begin{equation}
    \begin{aligned}
        \int_y^{2y} &\frac{\delta_p(s) \, ds}{s} = \int_y^{2y} \bigg[\frac{\C_p(\overline{B_s}\setminus \Omega;B_{2s)}}{\C_p(\overline{B_s}, B_{2s})} \bigg]^{\frac{1}{p-1}} \frac{ds}{s} \leq \gamma  \int_y^{2y} \bigg[\frac{\C_p(\overline{B_s}\setminus \Omega;B_{4y)}}{C_1 s^{N-p}} \bigg]^{\frac{1}{p-1}} \frac{ds}{s} \\
        & \leq ( 2^{N-p} \gamma/C_1)^{\frac{1}{p-1}}\int_y^{2y} \bigg[\frac{\C_p(\overline{B_{2y}}\setminus \Omega;B_{4y)}}{\C_p(\overline{B_{2r}}, B_{4r})} \bigg]^{\frac{1}{p-1}} \frac{ds}{s} = \gamma \delta_p(2y).      \end{aligned}
\end{equation} Hence for all $\N \ni m\ge i_l$ we have 
\[\int_{\sigma^{i_l} \rho_0}^{\rho_0} \frac{\delta_p(s) \, ds}{s}\leq \sum_{j=0}^{m-1} \int_{2^{-(j+1)}\rho_0}^{2^{-j} \rho_0} \frac{\delta_p(s)\, ds}{s} \leq \gamma \sum_{j=0}^{m-1}\delta_p(2^{-j} \rho_0).\]
\noindent The considerations done until this point are valid for all $p>1$, provided that condition \eqref{W} is satisfied with such exponent. \end{proof}

\subsection{Conclusion of the proof of Theorem \ref{boundary-estimate}}
\noindent The proof of Theorem \ref{boundary-estimate} hinges upon the possibility of finding a family of nested backward cylinders $\{Q_n\}_{n \in \N}$ centered at $(x_o,t_o) \in S_T$, where we can iteratively and quantitatively reduce the oscillation of the solution, truncated from above and below by the boundary datum. At this stage, the major difficulty of this double-phase parabolic problem is to deal with the method of the accommodation of its degeneracy. This is because of the double requirement due both to the intrinsic geometry and the restriction of the radii obliged by the phase, see Remark \ref{georeq} and condition \eqref{A}.

\subsection{Accommodation of the degeneracy}
\noindent Let $(x_o,t_o) \in S_T=\partial \Omega \times (0,T]$ and choose
\[k_0^+= \sup_{S_T} f\,, \quad \text{and} \quad k_0^-= \inf_{S_T} f\, , \]
\[\mu_0^{\pm}= \sup_{\Omega_T} (u-k^{\pm}_0)_{\pm}, \quad \text{and} \quad w_0^{\pm}=\mu_0^{\pm}- (u-k^{\pm}_{0})_{\pm}.\]
\[\omega_0= \osc_{\Omega_T} u\,.   \]
We assume $\mu_0^{\pm} > 0$, because otherwise there is nothing to prove. Now, for some $\epsilon\in (0,1)$ to be determined later, let us define for $s \in (0,1)$ the numbers
\begin{equation} \label{zerap}
    \tilde{\eta}_0(p,s)=  3\gamma^* [\delta_p(s)]^{{2-p}} s^{p-\epsilon} , \qquad \text{if}\qquad  \quad a(x_o,t_o)=0,
\end{equation}

\begin{equation}\label{zeraq}
         \tilde{\eta}_0(q,s)= 3 \gamma^* [\delta_q(s)]^{{2-q}} s^{q-\epsilon} , \qquad \text{if} \qquad  \quad a(x_o,t_o)>0,
\end{equation}

\noindent with $\gamma^*>0$  the geometric constant of Section \ref{geometric setting}, necessary for the application of Lemma \ref{lem3.3}.\vskip0.2cm 

\noindent Let us choose $\rho_0(p) \in (0,R_o)$ and $\rho_0(q) \in (0, \, \min \{R_o, R\}/24)$, with $R_o$ the number for which \eqref{A} is valid and $R$ the maximal radius \eqref{maximal-radius}, to be numbers that satisfy 
\begin{equation*}
    \begin{cases}\tilde{\eta_0} (p):= \tilde{\eta_0} (p,\rho_0(p))< \min \{t_o,\, R_o^2\},\\
    \mu_0^{\pm} \delta_p(\rho_0(p)) > 2 C_p \rho_0(p) \,
    \end{cases}
    \, \text{and} \quad \begin{cases}\tilde{\eta}_0 (q):=\tilde{\eta_0} (q, \rho_0(q))<\min \{t_o,\, R_o^2\},\\
\mu_0^{\pm} \delta_q(\rho_0(q)) > 4 C_q \rho_0(q)+ (4C_q)^{\frac{1}{q-1}} (\rho_0(q)/R),
    \end{cases}
\end{equation*}
\noindent where $C_p,C_q>0$ are the constants provided by Lemma \ref{lem3.3}. The existence of such $\rho_0(p), \rho_0(q)$ is guaranteed by assumptions \eqref{eq1.6}-\eqref{eq1.7} in each case.

\vskip0.2cm \noindent Indeed, let us show for instance  the case of positive phase: we suppose, by contradiction, that for all $s \in (0,\, \min \{ R_o, R\}/24)$ we have the alternative 
\[\tilde{\eta}_o(q,s) \ge \min \{t_o, R_o^2\} \quad \vee \quad \mu_0^{\pm} \delta_q(s) \leq 4 C_q s+ (4C_q)^{\frac{1}{q-1}} (s/R). \] 
Then, for every such $s$ we can estimate $\delta_p(s)$ from above
\[ \bigg(\frac{3\gamma^*}{\min\{t_o, R_o^2\}} \bigg)^{\frac{1}{q-2}} s^{\frac{q-\epsilon}{q-2}} +(4 C_q s+ (4C_q)^{\frac{1}{q-1}} (s/R))/\mu^{\pm}_0\ge [\delta_q (s) ]\,. \] 

\noindent Hence 
\[ \int_0^{R_o}  \delta_q(s) \frac{ds}{s} \leq \int_0^{R_o}\bigg(\frac{3\gamma^*}{\min \{t_o, R_o^2\}} \bigg)^{\frac{1}{q-2}} s^{\frac{2-\epsilon}{q-2}}ds+(4C_q/\mu_0^{\pm}) (1+1/R)\int_0^{R_o}\, ds <\infty,  \]
contradicting \eqref{eq1.7}.  
With such numbers $\{\tilde{\eta_0}(p), \rho_0(p) \}$ and $\{\tilde{\eta}_0(q), \rho(q)\}$ we define the cylinders
\begin{equation} \label{zeroP} \qquad Q_0(p, \pm)= B_{\rho_0(p)}(x_o) \times \bigg( t_o- [\mu_0^{\pm}]^{2-p}\tilde{\eta}_0(p)\, ,  \, t_o\bigg),\end{equation}

\begin{equation} \label{zeroQ} \qquad Q_0(q, \pm)= B_{\rho_0(q)}(x_o) \times \bigg( t_o- [\mu_0^{\pm}]^{2-q}\tilde{\eta}_0(q)\, ,  \, t_o\bigg).\end{equation}

\noindent From now on the proof is standard, we repeat it here for the sake of completeness, within a compact notation. \vskip0.1cm \noindent Let us indicate with an index $i \in \{p,q\}$ the radii $\rho_0(p), \rho_0(q)$ and the time-lengths $\tilde{\eta}_0(p), \tilde{\eta}_0(q)$ and the next quantities that we are going to define. Let 

\begin{equation*} \label{zeroA} \qquad Q_0(i)= B_{\rho_0(i)}(x_o) \times \bigg( t_o- \tilde{\eta}_0(i)\, , \,  \, t_o\bigg),\quad \qquad  \omega_0(i)= \osc_{Q_0(i)} u\, \, .\end{equation*}

\subsection{The Iteration, First Step}

For $i=p,q$, we choose levels
\[k_0^+(i)= \sup_{Q_{0}(i,+)\cap S_T} f\, \, , \quad \text{and} \quad k_0^-(i)= \inf_{Q_{0}(i,-)\cap S_T} f\, \, .\] 
\noindent We define
\[\mu_0^{\pm}(i)= \sup_{Q_0(i,\pm)} (u-k^{\pm})_{\pm}\, \, , \quad \text{and} \quad w_0^{\pm}(i)=\mu_0^{\pm}(i)- (u-k^{\pm}_0(i))\, .\]
\newline \noindent Now we observe that for both $i=p,q$ we can always assume \begin{equation}\label{cuggino}
(\mu_0^{\pm}(i))^{2-i} \, \, \rho_0^{i}\leq \rho_0^{i-\epsilon}.
\end{equation}
\noindent because otherwise the quantities $\mu_0^{\pm}(i)$ are smaller than a power of the radius $\rho_0(i)$, for $i=p,q$ respectively, and we are done. This means that  
\[Q_0(i, \pm) \subseteq B_{\rho_0(i)}(x_o) \times \bigg(t_o- \tilde{\eta}_0(i), \, t_o  \bigg),\qquad \qquad i=p,q,\]
and the special choice of $\rho_0(i)$ allows us to apply Lemma \ref{lem3.3} to get 
\[ \sup_{Q_{1,i}(\pm)} (u-k^{\pm}(i))_{\pm} \leq \mu_{1}^{\pm}(i), \]
for
\[\mu_1^{\pm}(i)= (1-1/(2C_i)) \mu_{0}^{\pm}(i), \]
\noindent where for $i=p,q$ we have defined
\begin{equation*}
    Q_{1,i}(\pm)= B_{\rho_0(i)/2}(x_o) \times \bigg(t_o- \gamma^* [\rho_0(i)]^i \bigg(\mu_0^{\pm}(i) \delta_i(\rho_0(i)) \bigg)^{2-i}/8, \, \, t_o \bigg)\quad.
\end{equation*}

\subsection {The Iteration, $n$-th Step}

\noindent Now, we consider now Lemma \ref{itera} with $\rho_o=\rho_0(i)$, $\mu_o=\mu_0^{\pm}(i)$, and $C_1=C_i$ for $i=p,q$ respectively, depending on the case the phase vanishes at $(x_o,t_o)$ or not.\vskip0.2cm 

\noindent With these stipulations, we can find two sequences of radii $\{\rho_{j,p}\}_{j \in \N}$, $\{\rho_{j,q}\}_{j \in \N}$ with $\rho_o(i)=\rho_0(i)$, satisfying to \eqref{1}-\eqref{2}-\eqref{3} and Remark \ref{rmkq}. We define 
\begin{equation}
    \begin{cases}
        \eta_{n,p}^{\pm}= \gamma^* \rho_{n,p}^p \bigg(\mu_0^{\pm} (p)\delta_p(\rho_{n,p}) \bigg)^{2-p}, & \text{if} \quad a(x_o,t_o)=0,\\
        \eta_{n,q}^{\pm} =\gamma^* \rho_{n,q}^q \bigg(\frac{\mu_0^{\pm}(q)}{a(x_o,t_o)}\delta_q(\rho_{n,q}) \bigg)^{2-q}, & \text{if} \quad a(x_o,t_o)>0,
    \end{cases}
\end{equation} \noindent and cylinders 
\begin{equation} Q_{n,i}(\pm)= B_{\rho_{n,i}}(x_o) \times (t_o-\eta_{n,i}^{\pm}\, ,\, t_o),\qquad \text{for} \quad i=p,q\, .\end{equation}
\noindent Let us suppose the assertion valid until step $(n-1)$ and let us prove it for step $n$.

\vskip0.3cm \noindent Within conditions \eqref{1}-\eqref{2} for $j=n-1$ we can apply Lemma \ref{lem3.3} and obtain, for  
\[\mu_n^{\pm}(i)= (1-1/(2C_i)) \mu_{n-1}^{\pm}(i), \qquad \text{for} \quad i=p,q,\] and
\[ \sup_{Q_{n,i}} (u-k^{\pm}(i))_{\pm} \leq  \mu_{n}^{\pm}(i), \]
\noindent where for $i=p,q$ we have defined
\[Q_{n,i}(\pm) =  B_{\rho_{n,i}}(x_o) \times \bigg(t_o- \gamma^* \rho_{n,i}^i \bigg(\mu_n^{\pm}(i) \delta_i(\rho_{n,i}) \bigg)^{2-i}/2, \, \, t_o \bigg),\]
and once observed that 
\[Q_{n,i}(\pm)\subseteq \,  B_{\rho_{n,i}/2}(x_o) \times \bigg(t_o-\gamma^* \rho_{n,i}^i [\mu_n^{\pm}(i)\delta_i(\rho_{n,i})]^{2-i}/8,\, t_o\bigg) \, .\]
Hence at the $n+1$-th step, the application of Lemma \ref{lem3.3} provides for $i=p,q$ the estimates
\begin{equation*}
    \begin{aligned}
        \sup_{Q_{n+1,i}(\pm)} (u-k^{\pm}(i))_{\pm} &\leq \mu_{n}^{\pm}(i)\,  \bigg(1-\frac{1}{C_i}\delta(\rho_{n,i})\bigg)\\
        &\leq \mu_0^{\pm}(i) \, \exp \bigg{\{} -\frac{1}{C_{i}} \sum_{j=1}^n \delta_i(\rho_{j,i})\bigg{\}}\leq \mu_0^{\pm}(i) \, \exp \bigg{\{} -\frac{1}{\gamma_4} \int_{\rho_{n+1,i}}^{\rho_0(i)} \delta_i(s)\frac{ds}{s}\bigg{\}},
    \end{aligned}
\end{equation*} \noindent with $\gamma_4=\gamma_4(C_{i})$, using Bernoulli's inequality
and \eqref{3}. Taking into consideration \eqref{cuggino}, this yields
\begin{equation}
     \sup_{Q_{n+1,i}(\pm)} (u-k^{\pm}(i))_{\pm} \leq \mu_0^{\pm}(i) \exp \bigg{\{} -\frac{1}{\gamma_3} \int_{\rho_{n+1,i}}^{\rho_0(i)} \delta_i (s) \frac{ds}{s} \bigg{\}}+ \gamma_3[\rho_0(i)]^{\frac{\epsilon}{q-2}},
\end{equation} \noindent and considering as usual for any $\rho\in (0,\, \,  \rho_0(i))$ an integer $n \ge 0$ such that $\rho_{n+1,i}\leq \rho \leq \rho_{n,i}$, we obtain 
\begin{equation}
     \sup_{Q_{\rho}(\mu_0^{\pm}(i))} (u-k^{\pm}(i))_{\pm} \leq \mu_0^{\pm}(i) \exp \bigg{\{} -\frac{1}{\gamma_3} \int_{\rho}^{\rho_0(i)} \delta_i (s) \frac{ds}{s} \bigg{\}}+ \gamma_3 [\rho_0(i)]^{\frac{\epsilon}{q-2}},
\end{equation} \noindent being for our choice of $n$,
\[Q_{\rho}(\omega_0) \subseteq Q_\rho(\mu_0^{\pm}(i)) = B_{\rho}(x_o)\times \bigg( t_o-(\mu_0^{\pm}(i))^{2-p} \rho^p\, , \, t_o \bigg) \subset Q_{n+1,i}(\pm),\]
where the first set inclusion is due to the degenerate exponent $q,p>2$ and the choice 
\[ \max\{\mu_{0}^+(i), \mu_{0}^-(i) \}\leq \omega_0(i) - \osc_{S_T\cap Q_0(i)} f \leq \omega_0. \]
\noindent Finally we combine the two aforementioned estimates for $k^{\pm}(i)$
to obtain 
\begin{equation}
     \osc_{Q_{\rho}(\omega_0)} u \leq \omega_0 \exp \bigg{\{} -\frac{1}{\gamma} \int_{\rho}^{\rho_0(i)} \delta_i (s) \frac{ds}{s} \bigg{\}}+ \osc_{S_T\cap Q_0(i)} f + 2\gamma_3 [\rho_0(i)]^{\frac{\epsilon}{i-2}},
\end{equation} \noindent and the proof is concluded.


\section*{Acknowledgements} \noindent S. Ciani acknowledges the support of the department of Mathematics of the University of Bologna Alma Mater, and of the PNR italian fundings 2021-2027. E. Henriques was financed by Portuguese Funds through FCT - Funda\c c\~ao para a Ci\^encia e a Tecnologia - within the Projects UIDB/00013/2020 and UIDP/00013/2020. I. Skrypnik is partially supported by the Simons Foundation (Award 1160640, Presidential Discretionary-Ukraine Support Grants, Skrypnik I.I.)
 
\vskip0.2cm 
\noindent {\bf Research Data Policy and Data Availability Statements.} All data generated or analysed during this study are included in this article.

\nocite{*}

\newpage

\section{Appendix} \label{Appendix}

\subsection{Proof of Lemma \ref{lem2.1}} 
\noindent Without loss of generality, let $(\bar{x}, \bar{t})$ be the origin in $\R^{N+1}$. We test \eqref{eq1.5} by $\zeta (u_{h}-k)_{-}$ and integrate over $(h, \tau)$, for $0<h< \tau <\eta-h$. Using conditions \eqref{eq1.2} and the continuity of $u$ as a map $[\tau,\eta] \rightarrow L^2(B_r)$, we let $h \downarrow 0$. This provides
\begin{equation}\label{ripiccio}
    \begin{aligned}
\int_h^{\tau} \int_{B_r} & \partial_t \zeta (u_h-k)_- \, dxdt\\
& \qquad \xrightarrow[h\downarrow 0] \qquad -\frac{1}{2} \int_{B_r} \zeta (u-k)_-^2\, dx \bigg|^{\tau}_0 + q \int_0^\tau \int_{B_r} (u-k)_-^2 ( \partial_t \zeta_2) (\zeta/ \zeta_2)\, dxdt:= \mathcal{I}_p ,\end{aligned}\end{equation} \noindent for all such $0<\tau<\eta$, and
\begin{equation*}
    \begin{aligned}
\int_h^{\tau} \int_{B_r} & [\mathbb{A}(x,t,\nabla u)]_h \bigg( \zeta \nabla (u_h-k)_-+ q(u_h-k)_- (\nabla \zeta) (\zeta/\zeta_1)\bigg) \, dxdt\\
& \qquad \xrightarrow[h\downarrow 0] \qquad \int_0^\tau \int_{B_r}  \mathbb{A}(x,t,u,\nabla u ) \bigg( -\zeta \nabla u + q(u-k)_- (\nabla \zeta) (\zeta/\zeta_1)\bigg) \chi_{[u<k]}\, \, dxdt:= \mathcal{I}_e,\end{aligned}\end{equation*}
with $ \mathcal{I}_p+\mathcal{I}_e \ge 0$. Manipulating the sign of this inequality together with the signs of its various terms, while using conditions \eqref{eq1.2}, we estimate the following energy term as 
\begin{equation*}
\begin{aligned}
\mathcal{I}=\sup\limits_{0<t< \eta}& \int_{B_{r}} \zeta (u-k)^{2}_{-}(x,t) \, dx +\gamma^{-1}\iint_{Q^{+}_{r,\eta}}\bigg( |\nabla [\zeta (u-k)_{-}]|^{p} + a(x,t) |\nabla [\zeta (u-k)_{-}]|^{q}\bigg) \,  dxdt
\\ &\leq \sup\limits_{0<t< \eta} \int_{B_{r}} \zeta (u-k)^{2}_{-}(x,t) \, dx +\gamma^{-1}\iint_{Q^{+}_{r,\eta}}\bigg( |\nabla (u-k)_{-}|^{p}\zeta \,+ a(x,t) |\nabla (u-k)_{-}|^{q}\zeta \bigg) dxdt + \\
&\qquad \qquad \qquad + \gamma^{-1}\iint_{Q^{+}_{r,\eta}}\varphi \bigg(x,t, |\nabla \zeta| (u-k)_- \bigg) \, dxdt=: E+\Phi
\\
&\leqslant
2q \int_0^\eta \int_{B_r} (u-k)_-^2 |\partial_t \zeta| \, dxdt+ \gamma \iint_{Q^{+}_{r,\eta}} (u-k)_- |\nabla \zeta| \bigg( |\nabla u|^{p-1} + a(x,t) |\nabla u |^{q-1}\bigg) \, dxdt+ \\
&\qquad \qquad \qquad + \iint_{Q^{+}_{r,\eta}}\varphi \bigg(x,t, |\nabla \zeta_1| (u-k)_- \bigg) \, dxdt.
\end{aligned}
\end{equation*} \noindent Using Young's inequality we notice that 
\[ \iint_{Q_{r,\eta}^+} (u-k)_- |\nabla \zeta| \bigg( |\nabla u|^{p-1} + a(x,t) |\nabla u |^{q-1}\bigg) \, dxdt\]
\[\leq  C(\epsilon) \Phi + \epsilon \iint_{Q_{r,\eta}^+}[ |\nabla (u-k)_-|^p\zeta+ a(x,t) |\nabla (u-k)_- |^q\zeta]  \, dxdt.\]
\noindent Hence, reabsorbing on the right-hand side the last terms, using the properties of $\zeta$ and the monotonicity of the function $\xi \rightarrow \phi(x,t, \xi)$ in the last variable, we get 
\[        \mathcal{I} \leq  \gamma \sigma^{-1} \frac{k^{2}}{\eta} |A^{-}_{k,r,\eta}|+\gamma \sigma^{-q}\iint_{A^{-}_{k,r,\eta}}\varphi(x,t,k/r)dxdt \leqslant
\gamma \sigma^{-q}\bigg(\frac{k^{2}}{\eta}+ [\varphi^{+}_{r,k}]\bigg)|A^{-}_{k,r,\eta}|. \]
In order to conclude, we observe that Young's inequality again can be used on the left-hand side as
\begin{equation*}
    \begin{aligned}
\bigg( 1&+a^- \bigg( \frac{k}{r} \bigg)^{q-p} \bigg) \iint_{Q_{r,\eta}^+} |\nabla [\zeta (u-k)_-]^p \, dxdt\\
&\leq  \iint_{Q_{r,\eta}^+} |\nabla [\zeta (u-k)_-]^p \, dxdt + \iint_{Q_{r,\eta}^+} a(x,t) \bigg( \frac{k}{r} \bigg)^{q-p} |\nabla [\zeta(u-k)_-]|^p \, dxdt \\
& \leq \iint_{Q_{r,\eta}^+} |\nabla [\zeta (u-k)_-]^p \, dxdt + 2^{-1}  a^+(x,t) \bigg( \frac{k}{r} \bigg)^{q} |A_{k,r,\eta}|^{-}+ 2^{-1} \iint_{Q_{r,\eta}^+} a(x,t) |\nabla [\zeta(u-k)_-]|^q\, dxdt\\
& \leq \gamma \bigg( \mathcal{I} + a^+(x,t) \bigg( \frac{k}{r} \bigg)^{q} |A_{k,r,\eta}|^{-} \bigg).
    \end{aligned}
\end{equation*}\noindent Inequality \eqref{eq2.2} centered at the origin is found by putting all the pieces of the puzzle together; then, the usual transformation of coordinates $y=\bar{x}+x$, $s=\bar{t}+t$ finishes the job.\newline 
The estimate \eqref{eq2.3} is proven similarly: choosing $\displaystyle (u_h-k)_- \zeta_1(x)$ as a test function, integral $\mathcal{I}_p$ gets simplified.

\subsection{Proof of Lemma \ref{lem2.2}}
\noindent Firstly, we assume $\delta>0$. We test \eqref{eq1.5} by $(u_{h}+\delta)^{-\alpha}\zeta$, $t\in (\bar{t}, \, \bar{t}+\tau-h)$, $0<\tau\leq \eta$, and integrate over $(\bar{t}, \, \bar{t}+\tau-h)$, 

\begin{equation}
    \begin{aligned}
0\leq I_{p,h}+ I_{e,h}=&\int_0^{\tau-h} \int_{B_r \times \{t\}}
 \bigg{\{} \partial_t u_h (u_h+\delta)^{-\alpha} \zeta +\\
 &+ [\mathbb{A}(x,u,\nabla u)]_h \bigg[-\alpha (u_h+\delta)^{-(1+\alpha)} (\nabla u_h) \zeta+ q (u_h+\delta)^{-\alpha} \zeta_2^q \zeta_1^{q-1} (\nabla \zeta_1)  \bigg]\bigg{\}}dx dt 
 \end{aligned}
\end{equation} \noindent 
Here we first notice that, by using Fubini-Tonelli theorem and chain rule for the weak time derivative we can rewrite $I_{p,h}$ as follows
\begin{equation}
    \begin{aligned}
        I_{p,h}= &\frac{1}{1-\alpha}\int_{B_r} \int_{o}^{\tau-h} \partial_t [(u_h+\delta)^{1-\alpha}\zeta]- (u_h+\delta)^{-\alpha} \partial_t \zeta \, dxdt\\
        &= \frac{1}{1-\alpha} \int_{B_r}[(u_h+\delta)^{1-\alpha}\zeta](x, t) \, dx \bigg|_{t=0}^{t=\tau-h}-\frac{1}{1-\alpha}\int_{B_r} \int_{o}^{\tau-h} (u_h+\delta)^{-\alpha}  \partial_t \zeta \, dxdt. 
    \end{aligned}
\end{equation} \noindent Now, in order  to let $h\downarrow 0$ we refer to the properties of Steklov approximation (see  for instance \cite{DiB}, Lemma 3.2 page 11): we use the fact that $u(t, \cdot) : [0, \eta] \rightarrow L^{1-\alpha}(B_r)$ is continuous, and use the structure conditions \eqref{eq1.2} with Young's inequality, in order to apply the dominated convergence theorem (and Fatou's one, on the left hand side) and get, by the generality of $0<\tau\leq \eta$, 
\begin{equation}
    \begin{aligned}
        &\sup_{\bar{t}< t< \bar{t}+\eta}\int_{B_r} [(u+\delta)^{1-\alpha}\zeta] (x,t) \, dx \leq\\
        &\frac{1}{1-\alpha} \int_0^{\eta} \int_{B_r} (u+\delta)^{1-\epsilon} \partial_t \zeta +\mathbb{A}(x,u,\nabla u ) \bigg[q (u+\delta)^{-\alpha} (\nabla \zeta) (\zeta/ \zeta_1)- \alpha (u+\delta)^{-(1+\alpha)} \zeta \nabla u \bigg] \, dxdt\\
        &  \\
        &  \leq \frac{\|\partial_t \zeta \|_{\infty}}{1-\alpha} \int_0^{\eta} \int_{B_r} (u+\delta)^{1-\epsilon}+ \int_0^{\eta} \int_{B_r} qB_2\bigg[(u+\delta)^{-\alpha} (\nabla \zeta) (\zeta/ \zeta_1)\bigg]\bigg( |\nabla u |^{p-1}+ a(x,t) |\nabla u |^{q-1}\bigg) dxdt \\
        &  \qquad \qquad \qquad \qquad  \qquad  \qquad \qquad \qquad-\alpha K_1  \int_0^{\eta} \int_{B_r} \bigg[(u+\delta)^{-(1+\alpha)} \zeta \bigg]\bigg( |\nabla u |^{p}+ a(x,t) |\nabla u |^{q}\bigg)\bigg{\}} \, dxdt
    \end{aligned}
\end{equation} \noindent 
In the second integral term, we use Young's inequality (weighted on $(u+\delta)^{-1-\alpha}\zeta_2^q$) for $(u+\delta) (\nabla \zeta)\zeta_1^{q-p}$  and $\epsilon |\nabla u|^{p-1} \zeta_1^{(p-1)}$ to conjugate powers $p$ and $p/(p-1)$, on the third integral term we use Young's inequality with the same weight for $(u+\delta) (\nabla \zeta)$  and $\epsilon |\nabla u|^{q-1} \zeta_1^{q-1}$ to conjugate powers $q$ and $q/(q-1)$. Choosing $\epsilon$ small enough to reabsorb these quantities on the fourth and fifth negative integral terms, we obtain  \begin{multline*}
\frac{1}{1-\alpha}\sup\limits_{\bar{t}<t<\bar{t}+ \eta}\int_{B_{r}(\bar{x})}(u+\delta)^{1-\alpha}\zeta\, dx +\frac{\alpha}{\gamma}\iint\limits_{Q^{+}_{r,\eta}(\bar{x}, \bar{t})}(u+\delta)^{-\alpha -1} |\nabla u |^p \zeta \,  dxdt + \\+ \frac{\alpha}{\gamma}\iint\limits_{Q^{+}_{r,\eta}(\bar{x}, \bar{t})}a(x,t)(u+\delta)^{-\alpha -1} |\nabla u|^{q}\zeta\,  dxdt \leqslant \frac{1}{(1-\alpha)}\|\partial_t \zeta \|_{\infty} \iint\limits_{Q^{+}_{r,\eta}(\bar{x}, \bar{t})}(u+\delta)^{1-\alpha}dxdt +\\ +\gamma \alpha^{1-p} \|\nabla \zeta \|_{\infty}^p \iint\limits_{Q^{+}_{r,\eta}(\bar{x}, \bar{t})}(u+\delta)^{p-\alpha -1}dxdt + \gamma \alpha^{1-q} \|\nabla \zeta \|_{\infty}^q a^{+}_{Q^{+}_{r,\eta}(\bar{x},\bar{t})}\iint\limits_{Q^{+}_{r,\eta}(\bar{x}, \bar{t})}(u+\delta)^{q-\alpha -1} dxdt.
\end{multline*} \noindent
The desired inequality is therefore obtained by noticing that 
\[\iint\limits_{Q^{+}_{r,\eta}(\bar{x}, \bar{t})} (u+\delta)^{-\alpha-1} |\nabla u |^{p} \zeta\, dxdt= \iint\limits_{Q^{+}_{r,\eta}(\bar{x}, \bar{t})} \bigg{\{} |\nabla [(u+\delta)^{\frac{p-\alpha-1}{p}} \zeta^{\frac{1}{p}}]|^p -\frac{1}{p} (u+\delta)^{p-\alpha-1} |\nabla \zeta|^p (\zeta/\zeta_1)^p\bigg{\}} \,dxdt,\]
and similarly with the third term on the left-hand side of \eqref{eq2.3}.\vskip0.1cm \noindent To include the case $\delta=0$, we let $\delta \downarrow 0$ in the obtained estimates \eqref{eq2.3}, concluding with the help of the Dominated Convergence Theorem.

\end{document}